\documentclass[12pt]{amsart}
\usepackage{amssymb}
\usepackage{latexsym}
\usepackage{epsf}
\newtheorem{thm}{Theorem}[section]
\newtheorem{Thm}{Theorem}

\newtheorem{lem}[thm]{Lemma}

\newtheorem{Cor}[Thm]{Corollary}
\newtheorem{Con}[Thm]{Conjecture}

\theoremstyle{definition}

\def\1{{\rm1\mathchoice{\kern-0.25em}{\kern-0.25em}
        {\kern-0.2em}{\kern-0.2em}I}}

\newcommand{\lmn}[1]{\vadjust{\setbox1=\vtop{\hsize 25mm
\parindent=0pt\baselineskip=9pt
\rightskip=4mm plus 4mm#1}
\hbox{\kern-26mm\smash{\raise .5ex\box1}}}}

\newcommand{\nc}{\newcommand}
\def\be#1\ee{\begin{equation}#1\end{equation}}
\nc{\bc}{\begin{center}} \nc{\ec}{\end{center}} \nc{\bb}{\mathbb}
\nc{\cal}{\mathcal} \nc{\frk}{\mathfrak} \nc{\N}{{\mathsf N}}
\nc{\K}{{\mathsf K}} \nc{\fk}{\mathbf{k}} \nc{\fn}{\mathbf{n}}
\nc{\fb}{\mathbf{b}}  \nc{\e}{\varepsilon} \nc{\ev}{{\rm{ev}}}

\theoremstyle{remark}
\newtheorem{rem}{Remark}

\def\Z{{\mathbb Z}}
\def\Q{{\mathbb Q}}
\def\N{{\mathbb N}}
\def\C{{\mathcal C}}
\def\s{{\mathfrak s}}
\def\R{{\mathbb R}}
\def\A{{\mathcal A}}
\def\v8{\vskip 8pt}
\def\a{\alpha}
\def\b{\beta}

\def\g{\gamma}
\def\G{\Gamma}

\def\O{\mathcal O}

\def\F{\mathcal F}

\newcommand{\Cob}{{\mathcal Cob}}

\newcommand{\Cobi}{{\mathcal Cob}_{/i}}
\newcommand{\Cobl}{{\mathcal Cob}_{/l}}

\newcommand{\CC}{{\bf\mathcal Cob}_f}
\newcommand{\CCi}{{\bf\mathcal Cob}_{f/i}}
\def\ll{\left[\!\!\left[}
\def\rr{\right]\!\!\right]}

\newcommand{\Ztwo}{{\Z_{(2)}}}
\newcommand{\Kob}{\operatorname{Kob}}
\newcommand{\Kobh}{{\operatorname{Kob}_{/h}}}
\newcommand{\Kobpr}{{\operatorname{Kob}_{/\pm h}}}
\newcommand{\Kom}{\operatorname{Kom}}
\newcommand{\Komh}{\operatorname{Kom}_{/h}}
\newcommand{\Mat}{\operatorname{Mat}}
\newcommand{\Kh}{\operatorname{Kh}}
\newcommand{\Kk}{{\bf \mathcal Kh}}
\newcommand{\Mor}{\operatorname{Mor}}
\newcommand{\Obj}{\operatorname{Obj}}

\usepackage{
  amsmath,graphicx,amssymb,color,hyphenat,
  multicol,slashbox,rotating,floatflt,
}
\usepackage[all]{xy}

\newcommand{\eps}[2]{{\hspace{-3pt}\begin{array}{c}%
  \raisebox{-2.5pt}{\includegraphics[width=#1]{#2.eps}}%
\end{array}\hspace{-3pt}}}
\newcommand{\epsg}[2]{{\hspace{-3pt}\begin{array}{c}%
  \raisebox{0pt}{\includegraphics[width=#1]{#2.eps}}%
\end{array}\hspace{-3pt}}}

\nc{\fs}{\mathbf{s}}

\begin{document}

\title[Categorification of the colored Jones
polynomial]{Categorification of the colored 
Jones polynomial\\[0.3cm] and Rasmussen invariant of  links}
\author{Anna Beliakova}
\email{anna@math.unizh.ch}
\author{Stephan Wehrli}
\email{stephan.wehrli@math.unizh.ch}
\address{Institut f\"ur Mathematik, Universit\"at Z\"urich,
 Winterthurerstrasse 190,
CH-8057 Z\"urich, Switzerland}

%\date{October 2005}
%\thanks{The first author is grateful}

\begin{abstract}
We define a family of formal Khovanov brackets
of a colored link depending on two parameters.
The isomorphism classes of these brackets are 
invariants of framed colored links.
The Bar--Natan functors applied to these brackets
produce Khovanov and Lee homology theories categorifying the colored
Jones polynomial. Further,
we study conditions under which
framed colored link cobordisms induce chain transformations between
our formal brackets. We conjecture  that,
for special choice of parameters,  Khovanov and Lee homology theories
of colored links are functorial (up to sign).
Finally, we extend the Rasmussen invariant to links and give examples,
where this invariant is a  stronger obstruction to sliceness
than the multivariable Levine--Tristram signature.
\v8\v8


\end{abstract}

\maketitle

%%%%%%%%%%%%%%%%%%%%%%%%%%%%%%%%%%%%%%%%%%%%%%%%%%%%%%%%%%%
%{\bf $\ddot{\smile}$ Corrections are printed in boldface and
%put between two smileys. $\ddot{\smile}$}
%%%%%%%%%%%%%%%%%%%%%%%%%%%%%%%%%%%%%%%%%%%%%%%%%%%%%%%%%%

\section*{Introduction}
In \cite{kh1},  Khovanov constructed a bigraded
chain complex, whose Euler characteristic is the Jones polynomial
and  whose
chain equivalence class is a link invariant. In particular,
the bigraded homology group, known as Khovanov homology,
is a link invariant.
 Bar--Natan \cite{BN1} and the second author \cite{W}
showed that Khovanov
homology is strictly stronger than the Jones polynomial.
Furthermore, Khovanov homology is functorial with respect to link
cobordisms smoothly embedded in $\R^4$.

In \cite{Lee}, Lee modified Khovanov construction and made it
more accessible for calculations. The generators of Lee homology are
known explicitly. The middle topological degree of the two
generators of Lee homology is a new  knot invariant
introduced by Rasmussen  \cite{Ra}.
Rasmussen used it to give a combinatorial proof of the Milnor conjecture.
Note that this conjecture was previously accessible only
 via gauge theory --
instanton Donaldson invariants,  Seiberg--Witten theory or
  Ozsv\'ath--Szab\'o knot  Floer homology.
Viewing Khovanov theory  as a combinatorial counterpart of
the knot Floer homology of Ozsv\'ath and Szab\'o,
one can expect that the categorification of quantum 3--manifold invariants
will provide a combinatorial approach to Heegaard Floer homology.

The first step in this direction is a categorification of the colored
Jones polynomial. In \cite{Kho}, Khovanov made two proposals
for such a homology theory, based on two natural normalizations
of the colored Jones polynomial. Unfortunately,
 the first homology theory categorifying the colored Jones
polynomial is defined over
$\Z/2\Z$ and the second one, for the reduced Jones polynomial,
works   for knots only.
In this paper, we develop both Lee and Khovanov homology theories
of colored links over $\Z[1/2]$. To do this, we explore the ideas of Bar--Natan
 \cite{BN}, who regard  these theories  as just different
functors applied to the {\it formal Khovanov bracket}.
A similar approach to constructing
new homology theories over $\Z/2\Z$ for colored links was independently
proposed by
Mackaay and Turner \cite{MT}.

Our main results are summarized in the next subsection.

\subsection{Main results}
Let $\fn=\{n_1, n_2,...,n_l\}$ be a finite sequence of natural numbers.
Let $L_\fn$ be an oriented framed colored  link of $l$  components, where
 $n_i$ is the color of the $i$--th component, and $D_\fn$ be its diagram
in blackboard framing.
%Let $J(L^\fn)$
%be the Jones polynomial of $\fn$--cable of $L$.

In Section 2
we define the formal Khovanov bracket
$\ll D_\fn\rr_{\a,\b}$ of the colored link $L_\fn$ as an
object of $ \Kom(\Mat(\Kobh))$.
Here $ \Kom(\Mat(\Kobh))$
is the category of formal complexes over a `matrix extension'
of the category $\Kobh$, where Bar--Natan's formal brackets of
links belong to (see Section \ref{bn}).
%%%%%%%%%%%%%%%%%%%%%%%%%%%%%%%%%%%%%%%%%%%%%%%%%%%%%%%%%%%%%%%%%%

$\Kobh$ is itself a homotopy category
of complexes, so we may think of
$\ll D_\fn\rr_{\a,\b}$ as a `complex of complexes'.
The subscripts $\alpha$ and $\beta$ are two integer parameters
which enter in the definition of the differential of
$\ll D_\fn\rr_{\a,\b}$.

%%%%%%%%%%%%%%%%%%%%%%%%%%%%%%%%%%%%%%%%%%%%%%%%%%%%%%%%%%%%%%%%%%
We show that

\begin{Thm}\label{one} For any $\a$ and $\b$,
the isomorphism class of the complex $\ll D_\fn\rr_{\a,\b}$ 
is an invariant of the colored framed oriented link $L_\fn$.
\end{Thm}

Let $\A$ be the category of $\Z[1/2]$--modules.
By applying the Khovanov functor $\F_{\Kh}$ 
and the Lee  functor $\F_{\rm Lee}$ 
to the formal bracket, we
get homology theories over $\Komh(\A)$. 

%%%%%%%%%%%%%%%%%%%%%%%%%%%%%%%%%%%%%%%%%%%%%%%%%%%%%%%%%%%%%%%%%%%%%%

\begin{Cor}
The total graded Euler characteristic of
$\F_{\Kh}(\ll D_\fn\rr_{1,0})$ is equal to
the colored Jones  polynomial of $L_\fn$.
\end{Cor}

The precise definition of the total graded Euler characteristic
will be given in Subsection \ref{firstapproach}.
% Let $\A$ be the category of $\Z/2\Z$--modules, then
% $\F_{\Kh}(\ll D_\fn\rr_{1,0})$ is the categorification of the colored
% Jones polynomial defined in \cite{Kho}.
Note that Khovanov's \cite{Kho} categorification of the colored
Jones polynomial is a variant of $\F_{\Kh}(\ll D_\fn\rr_{1,0})$
with coefficients in $\Z/2\Z$.

%In Section \ref{fr}
%we study the category of framed link cobordisms $\Cob_f^4$.
%The objects of this category are  diagrams of framed links, and the
%morphisms are 2--dimensional framed cobordisms between such diagrams,
%generically embedded in $\R^3\times [0,1]$.
%One can regard such cobordisms 
%as a  movie of link diagrams with marked points.
%We extend the Carter--Saito movie moves
%to this setting.

In Section \ref{fr} we study movie presentations
of framed cobordisms, where by a framed cobordism
we mean a compact smooth oriented surface which is properly embedded
in $\R^3\times I$ and equipped with a trivialization of its
normal bundle in $\R^3\times I$, and which connects a framed link
in $\R^3\times\{0\}$ to a framed link in $\R^3\times\{1\}$.
We extend the Carter--Saito movie moves \cite{CS}
to the setting of framed cobordisms.

%%%%%%%%%%%%%%%%%%%%%%%%%%%%%%%%%%%%%%%%%%%%%%%%%%%%%%%%%%%%%%%%%%%%%%

\begin{Thm}
Two movies present isotopic framed cobordisms
if and only if there is a sequence of modified 
Carter--Saito moves and additional moves depicted in Figure 4
that takes one movie to the other.
\end{Thm}

%Let $\CC^4$ be the category of colored framed cobordisms.
%The objects are now framed colored link diagrams and morphisms
%have to respect colors.

%%%%%%%%%%%%%%%%%%%%%%%%%%%%%%%%%%%%%%%%%%%%%%%%%%%%%%%%%%%%%%%%%%%%%%

A colored framed cobordism is a framed cobordism together
with a coloring of its connectivity components by natural numbers.
Colored framed cobordisms have movie presentations whose stills
are colored framed link diagrams. Let $\CC^4$ be the category
whose objects are colored framed link diagrams and whose morphisms are
movie presentations of colored framed cobordisms. In Sections 4 and 5
we show that $\F_{\Kh}(\ll D_n\rr_{0,1})$ and
$\F_{\rm Lee}(\ll D_n\rr_{1,1})$ extend to functors $\F_{\Kh}\circ\Kk_{0,1}$
and $\F_{\rm Lee}\circ\Kk_{1,1}$, respectively, from $\CC^4$ to the category
of complexes over $\A$. More precisely,

\begin{Thm} The functors
$\F_{\Kh}\circ\Kk_{0,1}$ and 
$\F_{\rm Lee}\circ\Kk_{1,1}$ from $\CC^4$ to $\Kom(\Komh(\A))$
are  well--defined.
\end{Thm}

%The proof is given in Section 5.
Let $\CCi^4$ be the quotient of $\CC^4$ 
by framed Carter--Saito movie movies, and $\Komh(\Komh(\A))_{/\pm}$
be the projectivization of $\Komh(\Komh(\A))$, where each morphism is 
identified with its negative.
 We expect 

\begin{Con} The functors $\F_{\Kh}\circ\Kk_{0,1}$ and 
$\F_{\rm Lee}\circ\Kk_{1,1}$ descend to functors
$\CCi^4\to \Komh(\Komh(\A))_{/\pm}$.
\end{Con}

%%%%%%%%%%%%%%%%%%%%%%%%%%%%%%%%%%%%%%%%%%%%%%%%%%%%%%%%%%%%%%%%%%%%%

Finally, we extend the definition of the Rasmussen invariant to links
and study its properties.
We  show that 
%this invariant provides a similar obstruction
%to sliceness as a multivariable Levine--Tristram signature. 
%It turns out, that
in some cases the Rasmussen invariant of links is a stronger 
obstruction to sliceness
than the multivariable  Levine--Tristram signature  defined 
 by Cimasoni and Florens \cite{CF}.
%%%%%%%%%%%%%%%%%%%%%%%%%%%%%%%%%%%%%%%%%%%%%%%%%%%%%%%%%%%%%%%%%%%%%

Another interesting application of the Rasmussen invariant of links
was found by Baader \cite{Ba}.
He used the Rasmussen invariant to define a quasimorphism on
the braid group and to estimate the torsion length for
alternating braids.

%%%%%%%%%%%%%%%%%%%%%%%%%%%%%%%%%%%%%%%%%%%%%%%%%%%%%%%%%%%%%%%%%%%%%

\subsection{Plan of the paper} 
In Section \ref{bn} we recall the Bar--Natan construction. Then
we define the formal Khovanov bracket of colored links.
In Section \ref{fr} we study framed cobordisms and their movie
presentations. Further, we construct  maps between  our formal
brackets induced by  colored framed link cobordisms.
The last section is devoted to the Rasmussen invariant of links.

\subsection*{Acknowledgment} 
 We wish to express our
gratitude to Dror Bar--Natan for allowing us to use his figures
from \cite{BN},
to Vincent Florens for teaching us the Levine--Tristram signature
 and
to Mikhail Khovanov for the e--mail conversation.

\section{Bar--Natan's construction}\label{bn}
In \cite{BN}, Bar--Natan defined
the formal Khovanov bracket $\ll\cdot\rr$ for any link (or tangle)
 in such a way  that the Khovanov and Lee's homology theories
can be reconstructed from $\ll\cdot\rr$. In this section
we briefly recall the Bar--Natan's construction.

\subsection{Formal Khovanov bracket}\label{11}
Suppose we have  a generic diagram $D$ of an oriented link $L$ in $S^3$
with $c$ crossings.
There is a  cube of resolutions associated to $D$ (compare \cite{BN1}).
The vertices of the cube
correspond to the configurations of circles obtained after smoothing
of all crossings in $D$.  For any crossing, two different smoothings are
allowed: the  $0$--  and  the $1$--smoothing. Therefore,
we have $2^c$ vertices. After numbering the crossings of $D$,
we can label the vertices of the cube by $c$--letter
strings of $0$'s and $1$'s,
specifying  the smoothing  chosen at each crossing.
The cube is skewered along its main diagonal,
from $00...0$ to $11...1$.
The number of 1 in the labeling of a vertex is
equal to its `height' $k$. The cube is displayed in
such a way that the vertices of height $k$ project down
to the point $r:=k-c_-$ (see Figure 1).

%\begin{figure}[p]
%\begin{center}
%\begin{sideways}
%  \input Main.pstex_t
 %\[ \includegraphics[width=6.4in]{figs/Main.eps} \]
%\end{sideways}
%\end{center}
%\end{figure}
\v8\v8

\begin{center}
%\hspace*{-2.6cm}
\mbox{\epsfysize=7cm \epsffile{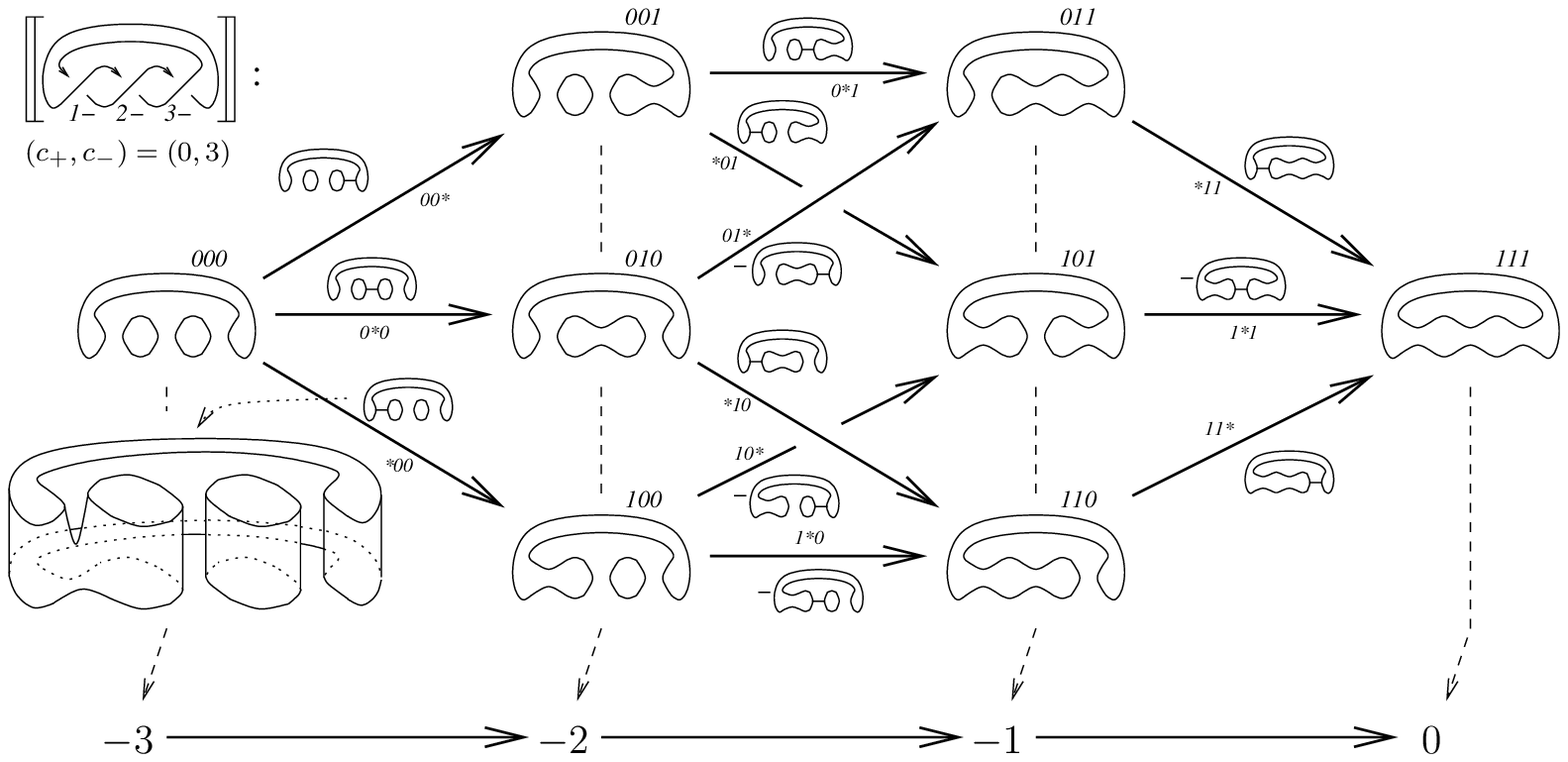}}
\v8

 Figure 1 {\it  The cube of resolutions  for the trefoil}
\end{center}
\v8

Two vertices of the cube are connected by an edge
if their labelings  differ by one letter. The edges are directed
(from the vertex where this letter is $0$ to the vertex where
it is $1$). The edges correspond to cobordisms from the tail
configuration of circles to the head configuration (compare Figure 1).

Bar--Natan proposed to interpret
the cube of resolutions denoted $\ll D\rr$
 as a complex, where
all smoothings are considered as spaces and all cobordisms as maps.
The $r$th chain space $\ll D\rr^{r}$ of the complex $\ll D\rr$
is a  formal direct sum of the
$\frac{c!}{k!(c-k)!}$ ``spaces'' at height $k$
in the cube and  the sum of ``maps''
with tails at height $k$
defines the  $r$th differential.

More precisely, $\ll D\rr$ is considered as an object of
$\Kom(\Mat(\Cob^3))$. Here $\Cob^3$ is the additive
category whose
objects are circle configurations (smoothings) and morphisms are
$2$--cobordisms between such smoothings.
For any additive category $\C$, $\Mat(\C)$ is the category
whose objects are formal direct sums of objects of $\C$ and whose
composition law is modeled on the matrix multiplication.
$\Kom(\C)$ is the category of complexes over $\C$, where objects are
chains of finite length and morphisms are chain transformations.

Let us impose  some local relations in $\Cob^3$:
(S) any cobordism  containing a closed sphere as a connected
component is set to be zero;
(T) any closed torus can be removed from a cobordisms at cost
of the factor $2$;
(4Tu) the four tube relation defined in \cite{BN}.
The neck cutting relation drawn below is a special form
of the 4Tu. If 2 is invertible, we can use this relation to cut
any tube inside a cobordism.
\v8

\begin{center}
\mbox{\epsfysize=1.2cm \epsffile{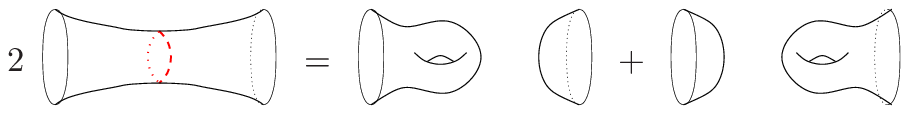}}
\v8

 Figure 2 {\it  The neck cutting relation}
\end{center}
\v8
We denote the quotient of $\Cob^3$ by these relations
$\Cobl^3$ and consider
$\ll D\rr$ as an object of  $\Kom(\Mat(\Cobl^3)$.
We set $\Kob:=\Kom(\Mat(\Cobl^3)$

\begin{thm}[Bar--Natan]\label{invar}
The homotopy type of $\ll D\rr $ is
an invariant of $L$.
\end{thm}

In \cite{BN}, Bar--Natan constructed explicit homotopies between complexes
 related by the three Reidemeister moves.

%%%%%%%%%%%%%%%%%%%%%%%%%%%%%%%%%%%%%%%%%%%%%%%%%%%%%%%%%%%%%%%%%%%%%%%%

\subsection{Topological grading}\label{grad}
A pre--additive category $\mathcal{C}$ is called
graded if it has the following additional properties.
Its morphism sets are graded Abelian groups, and the
degree is additive under composition of morphisms.
Moreover, there is a $\Z$--action $(m,\O)\mapsto\O\{m\}$
on the objects $\O$ of $\mathcal{C}$, which shifts the gradings
of the morphisms, but such that
$\Mor(\O_1\{m_1\},\O_2\{m_2\})=\Mor(\O_1,\O_2)$ as plain
Abelian groups.

Bar--Natan \cite{BN} observed that
$\Cobl^3$ can be transformed into a graded category
by introducing artificial objects $\mathcal{O}\{m\}$
for every $m\in\Z$ and every object $\mathcal{O}\in{\rm Obj}(\Cobl^3)$,
and by defining the degree of a cobordism $S\in\Mor(\Cobl^3)$
to be its Euler characteristic.
In what follows, we will denote by $\Cobl^3$ this
graded category, and we will refer to its grading as
the topological grading. Note that the topological
grading of $\Cobl^3$ induces topological
gradings on $\Mat(\Cobl^3)$ and $\Kom(\Mat(\Cobl^3))$.

%%%%%%%%%%%%%%%%%%%%%%%%%%%%%%%%%%%%%%%%%%%%%%%%%%%%%%%%%%%%%%%%%%%%%%%%

\subsection{Functoriality}\label{funct}

A link cobordism is 
is a compact oriented surface which is smoothly and properly
embedded in $\R^3\times I$ and connects a link in $\R^3\times\{0\}$
to a link in $\R^3\times\{1\}$. Splitting cobordisms into pieces
by planes $\R^3\times\{t\}$, $0\leq t\leq 1$, and projecting
down to the plane, we can view them as a sequence of link
diagrams or a movie of diagrams.
Altering $t$, we can  assume that
 any two consecutive diagrams in the movie differ by
 one of the following transformations
  --- a Reidemeister move, a cap or a cup, or a saddle.
It was shown in \cite{CS} that two such movies
present isotopic cobordisms if and only if they can be related
by a finite sequence of Carter--Saito movie moves.

Let $\Cob^4$ be the category whose objects are oriented
link diagrams, and whose morphisms are
movie presentations of cobordisms between links
described by such diagrams. Let $\Cobi^4$ be the quotient
of $\Cob^4$ by Carter--Saito movie moves.

The formal Khovanov bracket descends to a functor from
$\Kh:\Cob^4 \to \Kob$.
On the objects,  $\Kh(D)$ is defined as
the complex $\Kh^r(D):=\ll D\rr\{r+c_+-c_-\}$,
whose differentials are the same as those of $\ll D\rr$.
Note that all differentials in $\Kh(D)$
are of topological degree zero.
Moreover, it follows from the proof of Theorem \ref{invar} that
the graded homotopy type of $\Kh(D)$ is a link invariant
(cf. \cite[Theorem 3]{BN}).

On the generating morphisms $\Kh$  is defined as follows:
For the Reidemeister moves we take the chain homotopies constructed
 in \cite{BN} for the proof of Theorem \ref{invar}.
For the cup, cap or the saddle, we take the natural chain transformations
given by the corresponding cobordisms.

Let  $\Kobh$ be the category
 $\Kob$ modulo homotopies, i.e.  it has the same objects
as $\Kob$, but  homotopic morphisms
in $\Kob$ are identified.
Let $\Kobpr$ be the projectivization of $\Kobh$.

\begin{thm}[Bar--Natan]\label{two}
$\Kh$ descends to a functor $\Kh: \Cobi^4\to \Kobpr$.
\end{thm}

By the Carter--Saito theorem \cite{CS},  movie presentations of isotopic cobordisms
are related by  15 movie moves. Bar--Natan proved that
the  morphisms in $\Kob$ induced by these  movies moves   are
homotopic up to signs.

\subsection{Khovanov and Lee's theories}
Any functor from $\Cobl^3$ to an Abelian category
$\mathcal A$  extends to a functor $\mathcal F:\Kob\to \mathcal \Kom(\A)$
providing a homology theory. If in addition $\mathcal A$ is graded, and
$\mathcal F$ is degree--respecting, then the homology is a graded invariant
of a link.

\subsubsection{Khovanov functor}\label{Khfunctor}
Let $\O\in \Obj(\Cobl^3)$ and $\Ztwo=\Z[1/2]$.
We put
$$\mathcal F_{\Kh}(\O):=\Ztwo \otimes_{\Z} \Mor(\emptyset, \O)/ 
{\rm Rel}_{g>1}$$
where by ${\rm Rel}_{g>1}$ all
cobordisms of genus greater than 1 are set to be zero.
With a circle, $\mathcal F_{\Kh}$ associates the
$\Ztwo$--module of rank 2 generated by $v_+:=\eps{4mm}{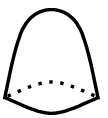}$
and by $v_-:=\frac12\eps{5mm}{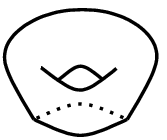}$.
The neck cutting relation allows to
identify the differentials in this theory
with the ones given by Khovanov \cite{kh1} (compare Exercise 9.3
in \cite{BN}). With the natural choice of 
grading on $\Ztwo$--modules  (${\rm deg}(v_+)=1$,  ${\rm deg}(v_-)=-1$),
the  functor $\F_{\Kh}$ is degree--respecting.
%%%%%%%%%%%%%%%%%%%%%%%%%%%%%%%%%%%%%%%%%%%%%%%%%%%%%%%%%%%%%%%%%%%%%%%%

Hence $\F_{\Kh}(\Kh(D))$ is a complex in the category of
graded $\Ztwo$--modules.
We define its graded Euler characteristic
$\chi(\F_{\Kh}(\Kh(D)))\in\Z[q,q^{-1}]$ by
$$
\chi(\F_{\Kh}(\Kh(D))):=\sum_{r,j} (-1)^rq^j{\rm
dim}_\Q(M^{r,j}(D)\otimes_{\Ztwo}\Q),
$$
where $M^{r,j}(D)$ denotes the homogeneous component of degree $j$
of the graded $\Ztwo$--module $\F_{\Kh}(\Kh^r(D))$.
It was shown in \cite{kh1} that $\chi(\F_{\Kh}(\Kh(D)))$ is equal to
the Jones polynomial of the link represented by the diagram $D$.

%%%%%%%%%%%%%%%%%%%%%%%%%%%%%%%%%%%%%%%%%%%%%%%%%%%%%%%%%%%%%%%%%%%%%%%%
\subsubsection{Lee's functor}
Let us put
$$\mathcal F_{\rm Lee}(\O):=\Ztwo \otimes_{\Z} \Mor(\emptyset, \O)/
(\epsg{9mm}{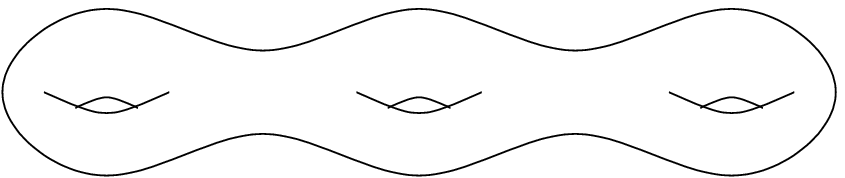}=8)$$
where the relation  set the morphism given by the genus 3
surface without boundary to be 8.
Here  the same rank 2 module  is associated to the circle.
But the differentials  $\Delta$ and $m$
are given by  the Lee's
formulas \cite{Lee}:
\be\label{diff}
  \Delta: \begin{cases}
    a \mapsto a\otimes a &\\
    b \mapsto b\otimes b &
  \end{cases}
  \qquad
  m_2: \begin{cases}
    a\otimes a\mapsto 2a &
    b\otimes b\mapsto -2b \\
    a\otimes b\mapsto 0 &
    b\otimes a\mapsto 0,
  \end{cases}
\ee
where  $a:=v_++v_-$ and $b:=v_+-v_-$.
The Lee functor is not degree--respecting.

\section{Formal Khovanov bracket of a colored link}
The aim of this section is to define the formal Khovanov bracket of a
 colored link. 
%categorifying the  colored Jones polynomial.
Our first approach is inspired by Khovanov \cite{Kho}.
Its modifications are necessary in order to get functoriality with respect to
colored framed link cobordisms.

\subsection{Colored Jones polynomial}
Let $\fn=\{n_1, n_2,...,n_l\}$ be a finite sequence of natural numbers.
Let $L_\fn$ be an oriented framed $l$  component link, whose
$i$--th component is colored by the $(n_i+1)$--dimensional
irreducible representation of $\frak sl_2$. Let $J(L^\fn)$
be the Jones polynomial of $\fn$--cable of $L$.
When forming the $m$--cable of a component $K$, we orient the
stands by alternating the original and the opposite directions.
More precisely,
let us enumerate the strands from left to right from 1 to $m$.
Then strand  1 is
oriented in the same way as $K$, the strand 2 is oppositely oriented, etc.

The colored Jones polynomial is given by the following formula.
\be\label{coljones}
J_\fn (L) =\sum^{\lfloor \fn/2\rfloor}_{\fk=\bf 0} (-1)^{|\fk|} \left(
\begin{array}{c} \fn-\fk \\ \fk\end{array}\right) J(L^{\fn-2\fk})
\ee
where $|\fk|=\sum_i  k_i$, and
$$\left(
\begin{array}{c} \fn-\fk \\ \fk\end{array}\right)=\prod^l_{i=1}
\left(
\begin{array}{c} n_i-k_i \\ k_i\end{array}\right)\, .$$

\begin{center}
\mbox{\epsfysize=10cm \epsffile{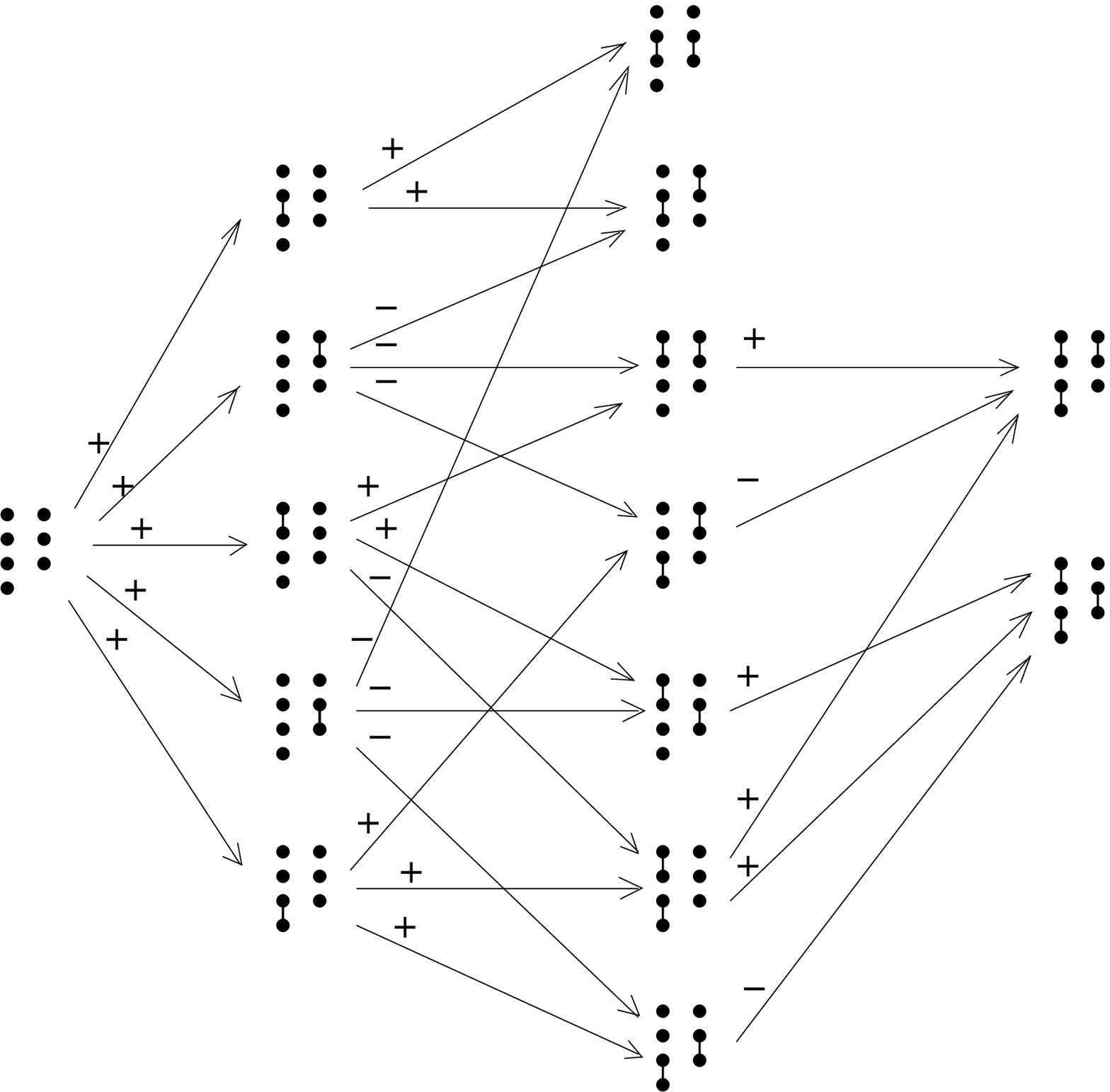}}
\v8

 Figure 3 {\it  The graph  $\G_{4,3}$.}
\end{center}

%\v8

\subsection{Graph $\G_\fn$}
The binomial coefficient
$\left(
\begin{array}{c} n-k \\ k\end{array}\right)$
equals the number of ways to select $k$ pairs of neighbors from $n$
dots placed on a line, such that each dot appears in at most one pair.
Analogously, $\left(
\begin{array}{c} \fn-\fk \\ \fk\end{array}\right)$
is the number of ways to select $\fk$ pairs of neighbors
on $l$ lines. We will call these choices $\fk$--pairings.

Let $\G_\fn$ be the graph, whose vertices correspond to $\fk$--pairings.
Two vertices of $\G_\fn$ are connected by an edge if the corresponding
pairings can be related to each other
 by adding/removing  one pair of neighboring points.
The height of a vertex labeled by a $\fk$--pairing is equal $|\fk|$.
The edges are directed towards increasing of heights (see Figure 3).

\subsection{Colored Khovanov bracket. First approach}\label{firstapproach}
Let $L_\fn$ be an oriented framed colored link as above and let $D_\fn$
be its generic diagram in blackboard framing.
Given $\G_\fn$ as above, we  associate to it
the formal Khovanov bracket $\ll D_\fn\rr$ of $L_\fn$
regarded as an element of $\Kom(\Mat(\Kobh))$. The construction
goes as follows.

At each vertex of $\G_\fn$
labeled by a $\fk$--pairing we put
the complex $\Kh(D^{\fn-2\fk})\in {\rm Obj}( \Kobh)$ defined in Subsection
\ref{funct}.
%%%%%%%%%%%%%%%%%%%%%%%%%%%%%%%%%%%%%%%%%%%%%%%%%%%%%%%%%%%%%%%%%%%%%%%

With an edge $e$ of $\G_\fn$ connecting $\fk$-- and $\fk'$--pairings,
we associate a morphism
$\Kh(e): \Kh(D^{\fn-2\fk})\to \Kh(D^{\fn-2\fk'})$
given by gluing of an annulus between the strands of the cable
which form a pair in $\fk'$, but not in $\fk$. According to the
definition of
$\G_\fn$ there is only one such pair.
Note that we view
the complexes $\Kh(D^{\fn-2\fk})$ and $\Kh(D^{\fn-2\fk'})$
as objects of the homotopy category $\Kobh$,
so that $\Kh(e)$ is a homotopy
class of chain transformations.
By Theorem \ref{two} this homotopy class is well--defined
up to sign. The sign of $\Kh(e)$ depends on the choice of the movie
presentation for the annulus. We call this choice satisfactory if
all squares of $\G_\fn$ anticommute.
Note that by Theorem \ref{two}
the squares of $\G_\fn$ commute up to sign,
because the cobordisms given by gluing of annuli in a different
order are isotopic.

Given a satisfactory choice of signs, the result is a complex
in $\Kom(\Mat(\Kobh))$, which we denote $\ll D_\fn\rr$.
The $i$--th chain of $\ll D_\fn\rr$
is a formal direct sum  of  complexes at height $i$,
i.e. $\ll D_\fn\rr^i:= \oplus_{|\fk|=i}\oplus_{\fs\in\fk}
\Kh(D^{\fn-2\fk})$,
where the notation $\fs\in\fk$ means that
$\fs$ is a $\fk$--pairing.
The $i$--th differential $d_i: \ll D_\fn\rr^i\to \ll D_\fn\rr^{i+1}$
is the formal sum of all morphisms $\Kh(e)$ corresponding to edges
with tails at height $i$.
%The category $\Kom(\Mat(\Kobh)) inherits the
%topological grading from $\Cobl^3$,
Because the Euler characteristic of an annulus is zero,
all $\Kh(e)$ have topological degree zero, and therefore
$\ll D_\fn \rr$ inherits a topological grading from
the topological gradings of the complexes $\Kh(D^{\fn-2\fk})$.
Besides the topological grading,
$\ll D_\fn \rr$ has two homological gradings, one
corresponding to the differential $d_i$ and one to the
differentials of the complexes $\Kh(D^{\fn-2\fk})$.
Note however that $\ll D_\fn \rr$ is not a bicomplex because
the chain transformations $\Kh(e)$ are considered up to homotopy.
It is an  interesting problem  whether one can  construct a bicomplex,
possibly by choosing suitable representatives for the
homotopy classes $\Kh(e)$. If such a bicomplex exists,
there should be a spectral sequence whose $E_2$ term
is determined by $\ll D_\fn \rr$ and which converges to the
homology of the total complex of that bicomplex.

We do not know how to form a bicomplex, but
we can define a total graded Euler characteristic
as follows.
Let $\ll D_\fn \rr^{i,r}\in{\rm Obj}(\Mat(\Cobl^3))$
be the formal direct sum 
$\ll D_\fn \rr^{i,r}:=\oplus_{|\fk|=i}\oplus_{\fs\in\fk} \Kh^r(D^{\fn-2\fk})$,
where $\Kh^r(D^{\fn-2\fk})$ denotes the $r$--th chain of
the complex $\Kh(D^{\fn-2\fk})$. The functor $\F_{\Kh}$
maps $\ll D_\fn \rr^{i,r}$ to a graded $\Z_{(2)}$--module
whose $j$--th homogeneous component we denote $M^{i,r,j}(D_\fn)$.
The total graded Euler characteristic of
$\F_{\Kh}(\ll D_\fn\rr)$ is defined by
$$
\chi(\F_{\Kh}(\ll D_\fn\rr)):=
\sum_{i,r,j}(-1)^{i+r} q^j{\rm dim}_\Q
(M^{i,r,j}(D_\fn)\otimes_{Z_{(2)}}\Q)\;.
$$
To complete this subsection, we prove the following lemma,
which shows that our construction of the colored Khovanov bracket
is well--defined.

%%%%%%%%%%%%%%%%%%%%%%%%%%%%%%%%%%%%%%%%%%%%%%%%%%%%%%%%%%%%%%%%%%%%%%%
%The functor $\F_{\Kh}$ applied to $\ll D_\fn\rr$ is degree--respecting.

\begin{lem}\label{signs}
For any graph $\G_\fn$ there exists a satisfactory choice of signs
making all squares anticommutative. Complexes defined
with different satisfactory sign choices are isomorphic.
\end{lem}

\begin{proof}
Let us first show that we can make all squares commutative.
We define a 1--cochain $\zeta\in C^1(\G_\fn, \Z/2\Z)$ as follows.
For any square $s\subset \G_\fn$,
we put  $\zeta(s)=1$ if $s$ is anticommutative and zero
otherwise.  We extend $\zeta$ by linearity to
$\G_\fn$. Now we multiply any map $\Kh(e)$
 by $(-1)^{\zeta(e)}$.

Note that $\zeta$ is well--defined, because there are no squares
which are  commutative and anticommutative simultaneously.
In other words,
the composition of maps induced by gluing of annuli is never zero.
Indeed, let $\kappa:\Kh(D^{\fn-2\fk})\rightarrow\Kh(D^{\fn-2\fk''})$
be a map induced by gluing $|\fk|-|\fk''|$ annuli. Let
$\bar{\kappa}:\Kh(D^{\fn-2\fk''})\rightarrow\Kh(D^{\fn-2\fk})$ denote
the map induced by the same annuli ``turned upside down''.
In the composition ${\kappa}\bar\kappa$, every annulus of
$\kappa$ is glued with the corresponding annulus of $\bar{\kappa}$,
such that the result is a torus. Hence ${\kappa}\bar{\kappa}$ is induced by
the union of $D^{\fn-2\fk''}\times[0,1]$ with a
collection of $|k|-|k''|$ tori. After isotopy, we can assume that these tori
lie in $\mathbb{R}^3\times \{1/2\}$. In $\mathbb{R}^3\times \{1/2\}$, the tori
may be linked with
$D^{\fn-2\fk''}\times\{1/2\}$, but if we consider
$2^{|\fk|-|\fk''|}{\kappa}\bar\kappa$ instead of ${\kappa}\bar\kappa$, we can
apply the neck cutting relation to obtain unlinked tori. It follows from
the (T) relation that $2^{|\fk|-|\fk''|}{\kappa}\bar\kappa$ is equal to
$4^{|\fk|-|\fk''|}$ times the identity morphism of
$\Kh(D^{\fn-2\fk''})$,
 and hence $\kappa$ is nonzero.

Given a complex with all squares commutative, we can make them
anticommutative  as follows. We multiply $\Kh(e)$
with $(-1)$ to the power number of pairings
to the right and above of the unique pairing in $\fk'\setminus\fk$.
 These signs are shown in Figure 3.

Given two satisfactory sign choices, the corresponding 1--cochains $\zeta$
and $\zeta'$ coincide on all squares, i.e. $\zeta-\zeta'=\delta\g$
with $\g\in C^0 (\G_\fn, \Z/2\Z)$. For any edge $e$ with boundary 
$s-s'$, we have $\zeta(e)- \zeta'(e)= \g(s)-\g(s')$. Therefore,
$(-1)^{\gamma}$ times  the identity map defines an isomorphism
between the corresponding complexes.\end{proof}

\begin{rem} Lemma \ref{signs} shows that the  categorification
of the colored Jones polynomial in \cite{Kho} can be defined over
integers.
\end{rem}

%\noindent{\bf Grading}

\subsection{Colored Khovanov bracket}\label{x}
In the following, we work with coefficients in
$\mathbb{Z}_{(2)}=\mathbb{Z}[1/2]$.
That is, we replace the category $\Cob^3$ of Section \ref{11}
by the category which has the same objects as $\Cob^3$ but
whose morphisms are formal $\mathbb{Z}_{(2)}$ linear combinations
of cobordisms.

Let us generalize the definition of $\ll D_\fn \rr$ 
 as follows.
As before,
we put  $\Kh(D^{\fn-2\fk})$
at  vertices  of $\G_\fn$
labeled by $\fk$--pairings. But
we modify the maps associated to edges of $\G_\fn$.
With an edge $e$ connecting $\fk$-- and $\fk'$--pairings
we associate the map $\Kh'(e):=\Kh(e)\circ(\a 1+\b X(e))$
where $1$ denotes the identity morphism of $\Kh(D^{\fn-2\fk})$ and
$X(e)$ is the endomorphism of $\Kh(D^{\fn-2\fk})$ defined below. Given a
satisfactory choice of signs, the result is a complex in $\Kom(\Mat(\Kobh))$,
which we denote $\ll D_\fn\rr_{\a,\b}$.
We have $\ll D_\fn\rr_{1,0}=\ll D_\fn\rr$. The functors
$\F_{\Kh}$ and $\F_{\rm Lee}$ can be applied to
$\ll \cdot\rr_{\a,\b}$ to obtain homology theories.
If $\b$ is nonzero, then
the topological degree of derivatives is not zero anymore.
%$\F_{\Kh}(\Kh'(e))$ does not respect the topological degree.
%We denote $\lr \cdot\rl:=\ll \cdot\rr_{1,1}$.

The map $X(e)$ is defined as follows.
 Assume $e$ is an edge between a $\fk$--pairing
and a $\fk'$--pairing, and let $C_i$ and $C_{i+1}$ be
the two strands of the cable of $L$ which form a pair in the
$\fk'$--pairing but not in the $\fk$--pairing. 
%Let $D^{\fn-2\fk}$ be a diagram of $L^{\fn-2\fk}$ and
Let us
 choose a point $P$ on $C_i$ which is not a crossing
of $D^{\fn-2\fk}$. Let $G$ be the region of $D^{\fn-2\fk}$
which lies next to $P$ and between the two components $C_i$ and
$C_{i+1}$. Color the regions of $D^{\fn-2\fk}$ in a chessboard
fashion, such that the unbounded region is colored white,
and put $\sigma(G):=+1$ if $G$ is black and
$\sigma(G):=-1$ if $G$ is white. Define
a cobordism $H(P)$ from $D^{\fn-2\fk}$ to itself as follows: $H(P)$
is the identity cobordism outside a small neighborhood of $P$, and
it is a composition of two saddle moves near $P$.
The first saddle splits off a small circle
from $C_i$. The second saddle merges the small
circle in $C_i$ again. Define $X(e):=(\sigma(G)/2)\Kh(H(P))$, where
$\Kh:\Cob^4\rightarrow\Kobh$ is the functor discussed in Section 
\ref{funct}.
We claim that $X(e)$ is
independent of the choice of the point $P$ on $C_i$. 
Indeed, moving the 
point $P$ past
a crossing of $D^{\fn-2\fk}$ changes the sign of both $\sigma(G)$ and
$\Kh(H(P))$. 
It is easy to see that $\F_{\rm Lee}(H(P))=-\F_{\rm Lee}(H(P'))$
if $P'$ is obtained from $P$ by moving past a crossings.
Moreover, the colors of regions next to $P$ and $P'$ are different.
If $C_i$ belongs to the cable of a component $K$ of $L$,
we also use the notation $X(K,i)$ for $X(e)$.

\subsection{\bf Proof of Theorem \ref{one}}
%%%%%%%%%%%%%%%%%%%%%%%%%%%%%%%%%%%%%%%%%%%%%%%%%%%%%%%%%%%%%%%%%%%%%

Let $D_\fn$ and $D'_\fn$ be two diagrams representing isotopic
colored framed links. Then $D^{\fn-2\fk}$ and $D'^{\fn-2\fk}$
represent isotopic links, and hence by
Theorem \ref{invar} the complexes
$\Kh(D^{\fn-2\fk})$ and $\Kh(D'^{\fn-2\fk})$
are isomorphic as objects of ${\rm Obj}( \Kobh)$.
The isotopy between the links represented by
$D^{\fn-2\fk}$ and $D'^{\fn-2\fk}$
extends to an isotopy between the annuli
appearing in the defintion of the differentials
of $\ll D_\fn\rr$ and $\ll D'_\fn\rr$.
Using Theorem \ref{two} and Lemma \ref{signs},
it easily follows that
$\ll D_\fn\rr$ and $\ll D'_\fn\rr$ are isomorphic.
$\hfill\Box$

\subsection{\bf Proof of Corollary 2}
The total graded Euler characteristic of $\F_{\Kh}(\ll D_\fn\rr)$ is
\begin{equation*}
\begin{split}
\chi(\F_{\Kh}(\ll D_\fn\rr))
&=\sum_{i,r,j}(-1)^{i+r}q^j{\rm dim}_\Q
(M^{i,r,j}(D_\fn)\otimes_{Z_{(2)}}\Q)\\
&=
\sum_i(-1)^i\sum_{|\fk|=i}\sum_{\fs\in\fk}
\chi(\F_{\Kh}(\Kh(D^{\fn-2\fk})))\\
&=\sum^{\lfloor\fn/2\rfloor}_{\fk=\bf 0} (-1)^{|\fk|}
 \left(
\begin{array}{c} \fn-\fk \\ \fk\end{array}\right)
\chi(\F_{\Kh}(\Kh(D^{\fn-2\fk})))\; .
\end{split}
\end{equation*}

Taking into account that
 $\chi(\F_{\Kh}(\Kh(D^{\fn-2\fk})))=J(L^{\fn-2\fk})$ we get the result.
$\hfill\Box$

%%%%%%%%%%%%%%%%%%%%%%%%%%%%%%%%%%%%%%%%%%%%%%%%%%%%%%%%%%%%%%%%%%%%%

\section{Framed cobordisms}\label{fr}

\subsection{Framings for submanifolds of codimension 2}
Let $M$ be a smooth oriented $n$--manifold and $N\subset M$ a compact smooth
oriented submanifold of $M$. By a framing of $N$ we mean a trivialization of
its normal bundle $\nu_N$ in $M$. Note that a smooth ambient isotopy between
submanifolds induces an isomorphism between their normal bundles. Hence it
makes sense to compare framings of ambient isotopic submanifolds. Given a
trivialization $f:\nu_N|_{\partial N}\rightarrow \partial N\times\mathbb{R}^2$,
we define a relative framing of $N$, relative to $f$, as a
trivialization of $\nu_N$ which restricts to $f$ on $\partial N$.
Relative isomorphism classes of oriented $2$--plane bundles over
$N$ which are trivialized over $\partial N$, correspond to homotopy classes of
maps from $(N,\partial N)$ to $(BSO(2),p_0)$, where $p_0$ is an arbitrary
basepoint in $BSO(2)$. Since $BSO(2)$
is a $K(\Z,2)$ space, we have $[N,\partial N;BSO(2),p_0]=H^2(N,\partial
N)=H_{n-4}(N)$. $N$ admits a relative framing if and only if $(\nu_N,f)$
corresponds to the zero class in $H_{n-4}(N)$. In that case, the set of all
relative framings is an affine space over $[N,\partial N;SO(2),1]=H^1(N,\partial
N)=H_{n-3}(N)$.

We are mainly interested in the case where $N$ is connected and $n=4$. In this
case the obstruction for the existence of relative framings is an integer
$e(\nu_N,f)\in H_0(N)=\mathbb{Z}$ which we call the relative Euler number of
$\nu_N$. The relative Euler number can be described explicitly as follows: let
$s$ be the zero section of $\nu_N$ and $s'$ a generic section such that
$f(s'(x))=(x,e_1)$ for $x\in \partial N$ where $e_1$ denotes the first basis
vector of $\mathbb{R}^2$. Then $e(\nu_N,f)=s\cdot s'$ where $s\cdot s'$ denotes
the algebraic intersection number of the surfaces $s$ and $s'$ in the total
space of $\nu_N$. $N$ has a tubular neighborhood in $M$ which is diffeomorphic
to the total space of $\nu_N$. Therefore, the relative Euler number $e(\nu_N,f)$
can be computed as a ``relative self--intersection number'' of $N$ in $M$.

\subsection{Framings for links and link cobordisms}\label{Flinks}

Let $K=N$ be a knot in $\R^3$. We can specify a framing of $K$ by
a vector field on $K$ which is nowhere tangent to $K$. If the vectors are
sufficiently short, their tips trace out a knot $K'$ parallel to $K$.
Recall that the framing coefficient $n(f)$ is defined as the linking number
of $K$ and $K'$.

Let us give an alternative description of the framing coefficient.
Let $S\subset\R^3\times I$ be a connected cobordism between the empty link and
the framed knot $K$, i.e. $\partial S=K\subset \R^3\times\{1\}$. We assume that
$S$ is parallel to the $I$ direction in a neighborhood of $\partial S$, such
that the restriction $\nu_S|_{\partial S}$ coincides with the normal bundle of
$\partial S$ in $\R^3\times\{1\}$. Then it makes sense to consider the relative
Euler number $e(\nu_S,f)$ where $f$ is the framing of $K$. We claim
that $e(\nu_S,f)=n(f)$.
We only prove that $e(\nu_S,f)$ is independent of the choice of $S$: let
$S,S_1$ be two cobordisms from the empty link to $K$ and let $\bar{S}_1$ denote
the cobordism $S_1$ ``turned upside down''. The composition of $S$ and
$\bar{S}_1$ is a closed surface $F:=S\cup \bar{S}_1$. Consider small
perturbations $S'$ and $S_1'$ of $S$ and $S_1$ with $\partial S'=\partial
S_1'=K'$ and let $F':=S'\cup\bar{S}_1'$. We have
$e(\nu_S,f)-e(\nu_{S_1},f)=S\cdot S' +\bar{S}_1\cdot\bar{S}_1'=F\cdot F'=0$,
where we have used that $F$ has self--intersection number zero because
$H_2(\R^3\times I)=0$. Hence $e(\nu_S,f)$ is independent of $S$.

Now let $S$ be a cobordism connecting two framed knots $(K_0,f_0)$ and
$(K_1,f_1)$. Choose cobordisms $S_0$ and $S_1$ from the empty link to $K_0$ and
$K_1$, respectively. By considering small perturbations $S',S_0',S_1'$ as above,
we obtain $0=S_0\cdot S_0'+S\cdot
S'+\bar{S}_1\cdot\bar{S}_1'=n(f_0)+e(\nu_S,f_0\cup f_1)-n(f_1)$. Hence $S$ 
admits a relative framing if and only if
$e(\nu_S,f_0\cup f_1)=0$ if and only if $n(f_0)=n(f_1)$.

Let us also consider the case of framed links. If $L$ is a link
of $|L|$ components in $\R^3$, a
framing of $L$ can be described by an $|L|$--tuple
$(n(f_1),\cdots,n(f_{|L|}))\in\Z^{|L|}$ where $f_i$ denotes the restriction of
the framing to the $i$th component. We define the total framing coefficient as
$$
n(f):=n(f_1)+\ldots +n(f_{|L|})+\sum_{i\neq j}lk(L_i,L_j).
$$
It is easy to see that $n(f)=e(\nu_S,f)$ for any connected cobordism $S$ from
the empty link to $L$. Arguing as above, we conclude that two framed links may
be connected by a relatively framed cobordism if and only if their total framing
coefficients agree.

If the set of relative framings
of $S$ is non--empty, it is an affine space over $H_1(S)$. The action of
$H_1(S)$ can be seen as follows: let $c$ be an oriented simple closed curve on
$S$ representing an element of $H_1(S)$. Consider a tubular neighborhood $U$ of
$c$, diffeomorphic to $c\times [0,2\pi]$. Let $\chi_c$ be the map from $S$ to
$SO(2)$ which is trivial on the complement of $U$ and maps a point
$(\theta,\varphi)\in U=c\times [0,2\pi]$ to rotation by $\varphi$. Then $c$ acts
on framings by sending the framing given by a vector field $v(z)$ to the framing
given by the vector field $\chi_c(z)v(z)$. In this context, the Poincar\'e dual
$\operatorname{PD}^{-1}[c]\in H^1(S,\partial S)$ has the following
interpretation: let $c'$ be a properly embedded
 simple curve on $S$ representing an element
of $H_1(S,\partial S)$. The restriction $\chi_c|_{c'}$ is a closed curve in
$SO(2)$ whose class in $\pi_1(SO(2),1)=\mathbb{Z}$ is given by
$[\chi_c|_{c'}]=c\cdot c'=\langle\operatorname{PD}^{-1}[c],[c']\rangle.$

\subsection{Link diagrams with marked points}

Let $L$ be a link and $D$ a diagram of $L$.
We may use $D$ to specify a framing on $L$,
namely the framing given by vector field on $L$
which is everywhere perpendicular to the plane of $D$.
This framing is called the blackboard framing.
It allows us to view link diagrams as diagrams of framed links.
The blackboard framing is invariant under the second and the third
Reidemeister moves, but not under the first Reidemeister move.
It is easy to see that two link diagrams describe isotopic framed links
if and only the differ by a sequence of the following moves:
the modified first Reidemeister move R1' shown below,
\begin{center}
\mbox{\epsfysize=3cm \epsffile{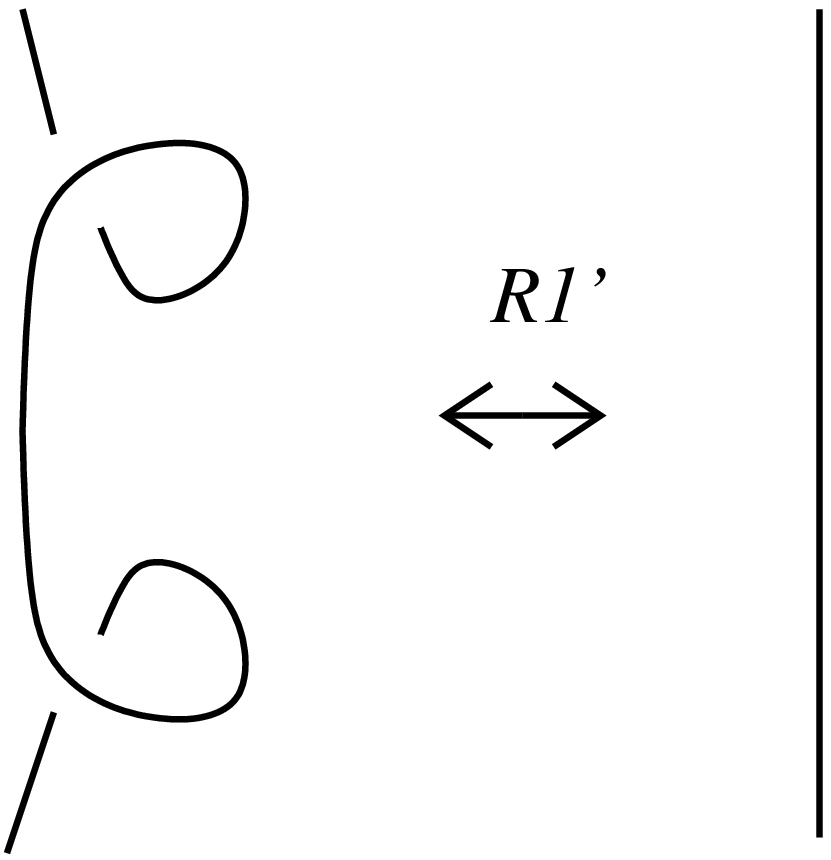}}
\end{center}
as well as the second and the third Reidemeister move.

A link diagram with marked points is a link diagram $D$ together with
 a finite
collection of distinct points, lying on the interiors of the edges of
$D$, and
labeled with $+$ or $-$. If $D$ is a link diagram with marked points, the
writhe $wr(D)$ is the difference between the numbers of positive and 
negative
crossings in $D$. The twist $tw(D)$ is the difference between the numbers of
positive and negative marked points. A link diagram with marked points
determines a framing $f_D$ of $L$ as follows: $f_D$ is given by a
vector field
which is perpendicular to the drawing plane, except in a small 
neighborhood of
the marked points, where it twists around the link, such that each positive
point contributes $+1$ to $n(f_D)$ and each negative point contributes 
$-1$ to
$n(f_D)$. Thus we have $n(f_D)=F\cdot L'=wr(D)+tw(D)$.

\begin{center}
\mbox{\epsfysize=2cm \epsffile{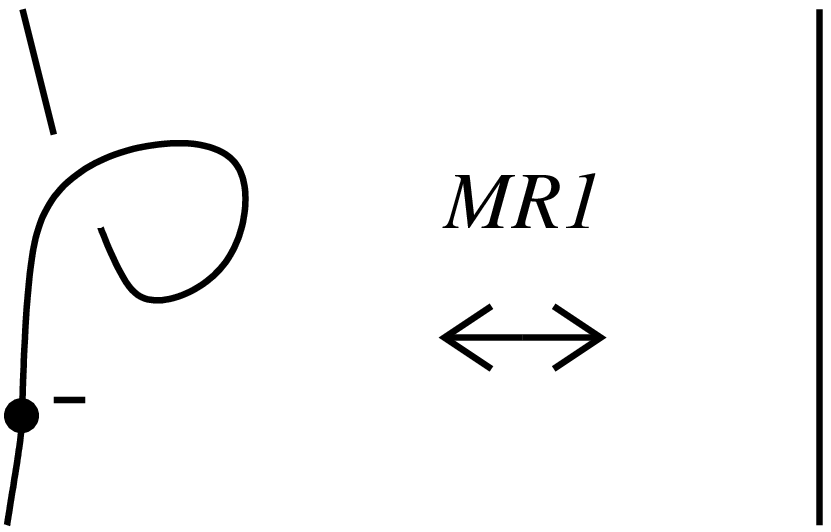}}
\end{center}

The marked first Reidemeister move MR1, shown above, leaves $n(f_D)$
unchanged. It
follows that two diagrams with marked points describe isotopic
framed links if
and only if they are related by a finite sequence of the
following moves: marked
first Reidemeister move MR1, Reidemeister moves R2 and R3,
creation/annihilation of
a pair of nearby oppositely marked points and sliding a marked point past a
crossing.

\v8
\begin{center}
%\hspace*{-2.6cm}
\mbox{\epsfysize=10cm \epsffile{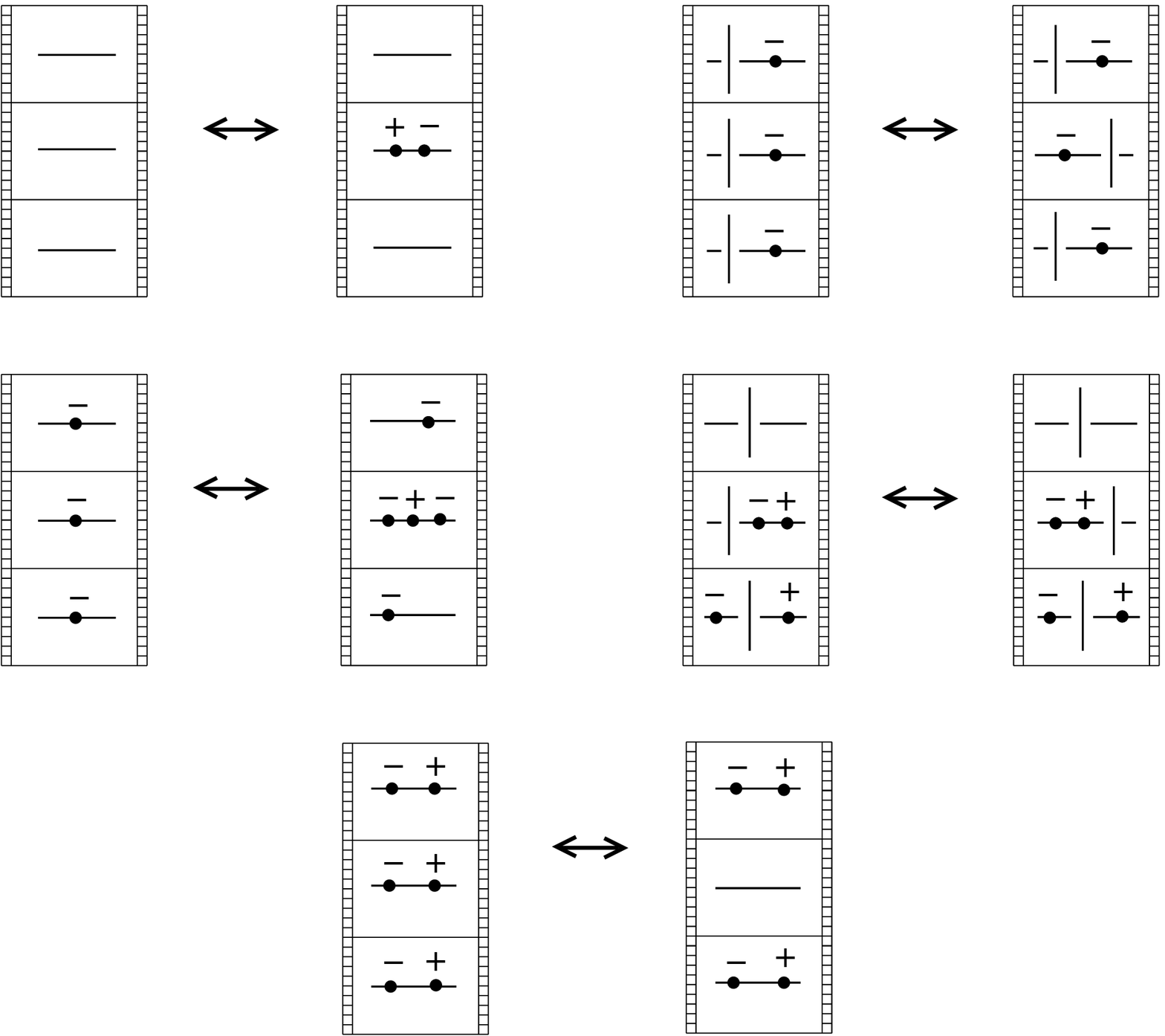}}
\v8 

 Figure 4 {\it  Additional movie moves for framed cobordisms}
\end{center}

\subsection{Movie presentations for framed cobordisms}

%Splitting cobordisms into pieces by planes $\R^3\times\{t\}$,
%$0\leq t\leq 1$, and projecting down to the plane, we can view
%them as a sequence of link diagrams or a movie of diagrams.
%We can assume that any two consecutive diagrams in the movie differ by
% one of the following transformations
%  --- a Reidemeister move, a cap or a cup, or a saddle.
In this subsection, we discuss movie presentations for framed link
cobordisms.

Let $S$ be an unframed cobordism, presented as a sequence of link diagrams.
If there are two consecutive link diagrams differing by an R1 move, we introduce
marked points in the movie presentation, such that every R1 move becomes
an MR1 move. The result is a movie of link diagrams with marked points,
describing a movie of framed links. The framings of these links determine a
well-defined framing of $S$. We claim that every framing of $S$ arises in this
way. To prove this claim, it would be sufficient to check that it is true
for elementary cobordisms (caps, cups and saddles).
However, we give a different proof. The marked points 
in the movie presentation
trace out curves on the cobordism $S$ (note that these curves may have local
extrema corresponding to annihilation and creation  of marked points).
We can orient these curves consistently, by declaring that positive points
``move'' in negative $I$ direction and negative points move in positive $I$
direction. Conversely, if $c$ is an oriented simple closed curve on
$S$, we can think of $c$ as consisting of lines traced out by marked
points. We can insert these marked points into a given movie presentation of
$S$. Thus, $c$ acts on movie presentations of $S$ by insertion of marked points.
This action induces an action on the framings of $S$ which are described
by the movie presentations. It is easy to see that the action on the framings
coincides with the
action of $H_1(S)$ discussed at the end of Subsection \ref{Flinks}. Since
$H_1(S)$ acts transitively, it follows that every framing of $S$ can be
described by a movie presentation. Moreover, two oriented simple closed curves
induce equivalent actions on framings if and only if they are homologous.

%\v8

\begin{proof}[\bf Proof of Theorem 3]
The local movie moves which are sufficient to relate any two
homologous curves traced out by marked points on a cobordism are
 shown in Figure 4. These movie moves,
together with modifications of the Carter-Saito movie moves obtained by
inserting marked points, are sufficient to relate any two 
movie presentations of isotopic framed cobordisms.
\end{proof}

We can transform a link diagram with marked points
into a link diagram without marked points
by inserting a left--twist curl for each point marked with a $+$ and a
right--twist curl for each point marked with a $-$. Under this substitution,
sliding a marked point past a crossing becomes a
composition of the second and the third Reidemeister move.
The MR1 move and creation/annihilation of a pair of oppositely marked
points become the modified first Reidemeister move R1'.

%Any framed cobordism can be viewed as
%a sequence of link diagrams without marked points, such that any two
%consecutive 
%diagrams differ by one of the moves just mentioned or by a cap, a cup or a
%saddle.

We can transform a movie presentation without marked points
into a movie presentation with marked points as follows: we replace each
R1' move by a composition of two opposite MR1 moves. The two opposite MR1 moves
create a pair of oppositely marked points. We annihilate this pair immediately
after its creation. The resulting movie presentation with marked points
differs from the original movie presentation only locally. 
Since we already know
movie moves for movie presentations with marked points, we can define
movie moves for movie presentations without marked points
simply by replacing marked points with curls.

\section{Colored framed  cobordisms}

Let $\CC^4$ be the category of colored framed movie presentations.
The objects are diagrams of colored links and the morphisms
movie presentations of colored framed links,
i.e. sequences of colored framed link diagrams,
where between two consecutive diagrams one of the following
transformations occur ---  R1', R2 or R3 move, a saddle, a cap or a cup.
Note that here we need to distinguish between two saddle moves:
a ``splitting'' saddle which splits one colored component into
two of the same color, and a ``merging'' saddle which
merges  two components of {\it  the same color} into one component.
To components colored differently the merging saddle can not be applied.
%Let $\CCi^4$ be $\CC^4$ modulo isotopies.

We are interested in a construction of a functor
%%%%%%%%%%%%%%%%%%%%%%%%%%%%%%%%%%%%%%%%%%%%%%%%%%%%%%%%
 $\Kk_{\alpha,\beta}:\CC^4\to \Kom(\Mat(\Kob_h))$.
%which descends to $\CCi^4$.
For the objects we put
 $\Kk_{\alpha,\beta}(D_\fn)=\ll D_\fn\rr_{\a,\b}$.
%%%%%%%%%%%%%%%%%%%%%%%%%%%%%%%%%%%%%%%%%%%%%%%%%%%%%%%%
The rest of this section is devoted to the definition of
chain transformations corresponding
to cap and cup, and saddles.
The  Reidemeister moves induce chain homotopies defined in \cite{BN}.

We introduce the following notation. Let $D_\fn$ be a
colored link diagram, and let $\fs$ be
a $\fk$--pairing of the $\fn$--cable of $D$.
Then we set $D^{\fs}:=D^{\fn-2\fk}$, where
the strands of $D^{\fn-2\fk}$ correspond to the dots which are
not contained in a pair of $\fs$.

\subsection{Cup and cap}\label{scupcap}
Consider two diagrams $D$ and $D_0$ which are related by a cap cobordism.
Assume that $D_0$ is the disjoint union of $D$ with a trivial component $K$.
Let $\fn_0$ be a coloring of $D_0$ and let $\fn$ denote the induced coloring
of $D$. Let $n$ denote the restriction of $\fn_0$ to $K$. Let
$\fs$ be a pairing of the $\fn$--cable of $D$ and let $\fs_0$ be a pairing
of the $\fn_0$--cable of $D_0$. We define a morphism
$\iota^{\fs_0,\fs}:\Kh(D^{\fs})\rightarrow\Kh(D_0^{\fs_0})$ as follows:
$\iota^{\fs_0,\fs}$ is nonzero only if the restriction of $\fs_0$ to $K$ is the
empty pairing (no pairs) and if $\fs_0$ agrees with $\fs$ on all other
components. In this case, we define $\iota^{\fs_0,\fs}$ as the composition of
the following two morphisms: the morphism induced by a union of $n$ caps whose
boundaries are the $n$ strands of the $n$--cable of $K$, and the endomorphism
$\varphi$ of $\Kh(D_0^{\fs_0})$ given by $\varphi:=\sum_{j=1}^n A_j\circ B_j$,
where $A_j$ denotes the composition of all morphisms $(\a 1-\b X(K,i))/2$ for
$1\leq i\leq j$, and $B_j$ denotes the composition of all morphisms $(\a 1+\b
X(K,i))/2$ for $j<i\leq n$.

Now let $D$ and $D_0$ be two link diagrams related by a cup cobordism.
Assume that $D$ is the disjoint union of $D_0$ with a trivial
component $K$. Let $\fn$ be a coloring of $D$ and let $\fn_0$ denote
the induced coloring of $D_0$. Let $n$ denote the restriction of $\fn$
to $K$. Let $\fs$ be a pairing of the $\fn$--cable of $D$ and
let $\fs_0$ be a pairing of the $\fn_0$--cable of $D_0$. We define a morphism
$\epsilon^{\fs_0,\fs}:\Kh(D^{\fs})\rightarrow\Kh(D_0^{\fs_0})$ as follows:
$\epsilon^{\fs_0,\fs}$ is nonzero only if the restriction of $\fs$ to $K$ is
the empty pairing and if $\fs$ agrees with $\fs_0$ on all other
components. In this case, we define $\epsilon^{\fs_0,\fs}$ as the composition
of the following two morphisms: the endomorphism $\varphi$ defined as above and
the morphism induced by $n$ cups whose boundaries are the $n$ strands of the
$n$--cable of $K$.

\subsection{Merging saddle}\label{smsaddle}

 Consider two diagrams $D$ and $D_0$
which are related by a saddle merging two components
$K_1$ and $K_2$ of $D$ into a single component $K$ of $D_0$.
Let $\fn$ be a coloring of $D$, such that
$K_1$ and $K_2$ have the same color $n$. Let $\fn_0$ be
the induced coloring of $D_0$. Consider a pairing $\fs$ of
the $\fn$--cable of $D$, and let $s_1$ and $s_2$
denote the restrictions of $\fs$ to $K_1$
and $K_2$, respectively. Let $s_1s_2$ denote the union of $s_1$ and
$s_2$, i.e. the pairing which consists of all pairs which are
contained in either $s_1$ or $s_2$.

Let $\g,\delta\in\mathbb{Z}$. Given a pairing $\fs_0$ of the
$\fn_0$--cable of $D_0$, the morphism
$$\psi^{\fs_0,\fs}_{\g,\delta}:\Kh(D^{\fs})\longrightarrow \Kh(D^{\fs_0})$$
 is nonzero only if the following is
satisfied: 
\begin{itemize}
\item
$s_1$ and $s_2$ have no common dot (meaning that there is no dot
which belongs to a pair both in $s_1$ and in $s_2$), 
\item
 $\fs_0$ is the pairing
which restricts to $s_1s_2$ on $K$ and which agrees with $\fs$ on all other
components. 
\end{itemize}
In this case, we define $\psi^{\fs_0,\fs}_{\g,\delta}$
as follows. For each pair of
$s_2$, consisting of dots numbered $i$ and $i+1$, consider the endomorphism $(\g
1+\delta X(K_1,i))/2$ of $\Kh(D^{\fs})$. Similarly, for each pair of $s_1$,
consisting of dots numbered $i$ and $i+1$, consider the
endomorphism $(\g 1+\delta X(K_2,i))/2$ of $\Kh(D^{\fs})$.
Denote the composition of these endomorphisms by $\psi_1$. Let $\fs'$ be the
$\fn$--pairing which restricts to  $s_1s_2$ on both
$K_1$ and $K_2$ and which agrees with the $\fn$--pairing $\fs$ on all other
components of $D$. Define a morphism $\psi_2$ from $\Kh(D^{\fs})$ to
$\Kh(D^{\fs'})$ as follows: for each pair of $s_2$, consisting of
dots numbered $i$ and $i+1$, consider an annulus attached
to the strands numbered $i$ and $i+1$ of $K^{s_1}$. Similarly,
for each pair of $s_1$, consisting of dots numbered $i$ and $i+1$,
consider an annulus attached to the strands numbered $i$ and $i+1$
of $K^{s_2}$. Let $\psi_2$ be the morphism induced by these annuli.
For every strand of $K_1^{s_1s_2}$ there is a corresponding
strand in $K_2^{s_1s_2}$. Let $\psi_3$ be the morphism from $\Kh(D^{\fs'})$
to $\Kh(D_0^{\fs_0})$ induced by merging each pair of corresponding strands
on $K_1^{s_1s_2}$ and $K_2^{s_1s_2}$ by a saddle cobordism. Define
$\psi^{\fs_0,\fs}_{\g,\delta}$ as the composition
$\psi^{\fs_0,\fs}_{\g,\delta}:=\psi_3\psi_2\psi_1$.

Note that our definition of the morphism $\psi^{\fs_0,\fs}_{\g,\delta}$
mimics the definition of the map $\psi$ in \cite{Kho}. 
Khovanov's map $\psi$
corresponds to our morphism $\psi^{\fs_0,\fs}_{0,2}$,
with the difference that we work with the Khovanov bracket
whereas Khovanov worked with Khovanov homology over
$\Z/2\Z$ coefficients.
Note that $\F_{\Kh}(\psi^{\fs_0,\fs}_{0,\delta})$ is
graded of degree $-n$ (where $n$ is the color of the
components involved in the saddle move, see above).
We denote by $\psi_{\g,\delta}$ the collection
of all morphisms $\psi^{\fs_0,\fs}_{\g,\delta}$.

\subsection{Splitting saddle}\label{sssaddle}
Suppose the diagrams $D$ and $D_0$ are related by a
saddle which splits a component $K$ of $D$ into two
components $K_1$ and $K_2$ of $D_0$. Let $\fn$ be
a coloring of $D$, and let $\fn_0$ be the induced
coloring of $D_0$. Consider a pairing $\fs$ of
the $\fn$--cable of $D$
which restricts to a $k$--pairing $s$ on $K$.

Let $\g,\delta\in\mathbb{Z}$.
Given a pairing $\fs_0$ of the $\fn_0$--cable of $D_0$,
the morphism
$$\bar{\psi}^{\fs_0,\fs}_{\g,\delta}:\Kh(D^{\fs})
\longrightarrow \Kh(D_0^{\fs_0})$$
is zero
unless $\fs_0$ has the following properties:
\begin{itemize}
\item
 the restrictions $s_1$ and $s_2$ of
$\fs_0$ to $K_1$ and $K_2$ have no common dot,
\item
 the union of $s_1$ and $s_2$ is equal to $s$, 
\item
 $\fs_0$ agrees with $\fs$ on all components of $D_0$ other
than $K_1$ and $K_2$. 
\end{itemize}
For an $\fs_0$ with these properties,
we define
$\bar{\psi}^{\fs_0,\fs}_{\g,\delta}:=2^k\bar{\psi}_1\bar{\psi}_2\bar{\psi}_3$,
where $\bar{\psi}_1,\bar{\psi}_2$ and $\bar{\psi}_3$ are the morphisms obtained
by turning the morphisms $\psi_1,\psi_2$ and $\psi_3$ of Subsection
\ref{smsaddle} upside down.

Note that $\bar{\psi}^{\fs_0,\fs}_{0,\delta}$ is graded of
degree $-n$, where $n$ is the color of the components involved
in the saddle move. We denote by $\bar{\psi}_{\g,\delta}$
the collection of all morphisms
$\bar{\psi}^{\fs_0,\fs}_{\g,\delta}$.

\subsection{Chain transformations induced by saddles}

Consider two diagrams $D$ and $D_0$ which are related by a merging
saddle. Assume we are given
a collection of morphisms $\psi=\{\psi^{\fs_0,\fs}\}$ as in Subsection
\ref{smsaddle} (we drop the subscripts $\g,\delta$ to simplify
the notation). We wish to have a criterion under which $\psi$ induces
a chain transformation from $\ll D_\fn\rr_{\a,\b}$ to $\ll
D_{0,\fn_0}\rr_{\a,\b}$.

Let $d$ and $d_0$ denote the
differentials of $\ll D_\fn\rr_{\a,\b}$ and
$\ll D_{0,\fn_0}\rr_{\a,\b}$, respectively.
Both $d_0\psi$ and $\psi d$
increase the height (the homological degree) by one.
Let $\fs$ be a pairing of the $\fn$--cable of $D$, and
let $\fs_0'$ be a pairing of the $\fn_0$--cable of $D_0$
whose height is one larger than the height of $\fs$.
Let $(d_0\psi)^{\fs_0',\fs}$ denote the ``restriction''
of $d_0\psi$ to $\Kh(D^\fs)$ and $\Kh(D_0^{\fs_0'})$,
and let $(\psi d)^{\fs_0',\fs}$ denote the ``restriction'' of
$\psi d$ to $\Kh(D^\fs)$ and $\Kh(D_0^{\fs_0'})$.
Assume that at least one of the morphisms $(d_0\psi)^{\fs_0',\fs}$
and $(\psi d)^{\fs_0',\fs}$ is nonzero. This is only
possible if the restrictions of $\fs$ to the two
components involved in the saddle move have no common dot.
Moreover, all pairs of the pairing $\fs_0$
(defined as in Subsection \ref{smsaddle}) must also
be pairs of $\fs_0'$. Therefore, $\fs_0'$ must contain a
unique pair $\pi$ which is not contained in $\fs_0$.
We assume that $\pi$ belongs to the
component of $D_0$ which is involved in the saddle move
(for otherwise $(d_0\psi)^{\fs_0',\fs}=\pm(\psi d)^{\fs_0',\fs}$
is trivially satisfied). Then we are in the
situation of (\ref{emerge}), where the pair in the
lower right corner is the pair $\pi$, and where we have
left away all dots corresponding to strands on
which $(d_0\psi)^{\fs_0',\fs}$ and $(\psi d)^{\fs_0',\fs}$
agree already by definition.

\begin{equation}\label{emerge}
\begin{split}
\epsfysize=3.8cm \epsffile{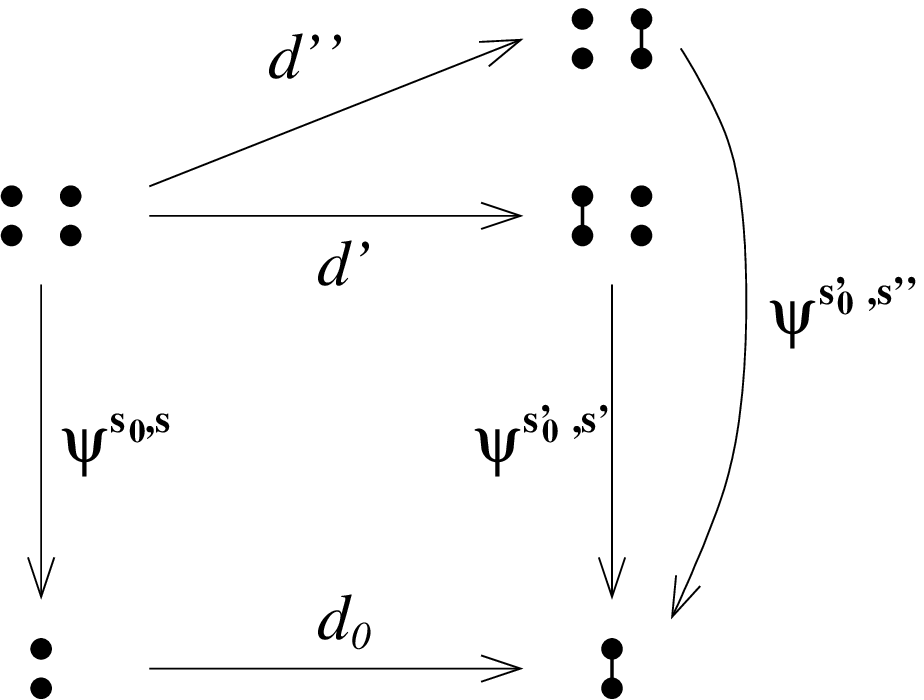}
\end{split}
\end{equation}

\v8

\begin{lem}\label{lmsaddle}
Assume that $d_0\psi^{\fs_0,\fs}=\pm
(\psi^{\fs_0',\fs'}d'+\psi^{\fs_0',\fs''}d'')$ for all
squares as in (\ref{emerge}). Then there is a $0$--cochain $\g\in
C^0(\G_{\fn},\mathbb{Z}/2\mathbb{Z})$ such that the morphisms
$(-1)^{\g(\fs)}\psi^{\fs_0,\fs}$ determine a chain transformation from
$\ll D_\fn \rr_{\a,\b}$ to $\ll D_{0,\fn_0}\rr_{\a,\b}$.
\end{lem}
\begin{proof}
Consider the subgraph $\G'_{\fn}$ of $\G_{\fn}$ whose
vertices are precisely those $\fn$--pairings $\fs$ 
whose restrictions to the two
components involved in the saddle move have no common dot. Let
$f:\G'_{\fn}\rightarrow \G_{\fn_0}$ be the map which maps a 
pairing/vertex $\fs$
to the induced pairing/vertex $\fs_0$. Note that every edge $e_0$ of
$\G_{\fn_0}$ appears as the lower edge of a square as in (\ref{emerge}). 
Define
a $1$--cochain $\zeta\in C^1(\G_{\fn_0},\mathbb{Z}/2\mathbb{Z})$ by setting
$\zeta(e_0):=0$ if $d_0\psi^{\fs_0,\fs}=+
(\psi^{\fs_0',\fs'}d'+\psi^{\fs_0',\fs''}d'')$ and $\zeta(e_0):=1$ if
$d_0\psi^{\fs_0,\fs}=-
(\psi^{\fs_0',\fs'}d'+\psi^{\fs_0',\fs''}d'')$. Since $f$ maps squares
of $\G'_{\fn}$ to squares of $\G_{\fn_0}$, and since all squares of $\G'_{\fn}$
and $\G_{\fn_0}$ anticommute, the $1$--cochain $f^*\zeta\in
C^1(\G'_{\fn},\mathbb{Z}/2\mathbb{Z})$ maps all squares of $\G'_{\fn}$ to zero.
Hence $f^*\zeta=\delta\g'$ for a $0$--cochain $\g'\in
C^0(\G'_{\fn},\mathbb{Z}/2\mathbb{Z})$. Now the $0$--cochain
$\gamma\in C^0(\G_{\fn},\mathbb{Z}/2\mathbb{Z})$ in the statement of the lemma
is an arbitrary extension of $\gamma'$. \end{proof}

Now assume $D$ and $D_0$ are related by a splitting
saddle and assume we are given a
collection of morphisms  $\bar{\psi}=\{\bar{\psi}^{\fs_0,\fs}\}$
as in Subsection \ref{sssaddle}.
Let $\fs$ be a pairing
of the $\fn$--cable of $D$ and let $\fs_0'$ be a pairing of the
$\fn_0$--cable of $D_0$, such that at least one of the morphisms
$(d_0\bar{\psi})^{\fs_0',\fs}$ and $(\bar{\psi} d)^{\fs_0',\fs}$
is nonzero. Let $K$ be the component of $D$ which is involved
in the saddle and let $s$ be the restriction of $\fs$ to $K$.
Similarly, let $K_1$ and $K_2$ be the components of $D_0$ which
are involved in the saddle and let $s_1'$ and $s_2'$ denote the
restrictions of $\fs_0'$ to $K_1$ and $K_2$.
Then every pair of $s$ must also appear in the
union $s_1'\cup s_2'$ (where we regard $s_1'$
and $s_2'$ as sets of pairs). If $s_1'$
and $s_2'$ have a common pair, we are in the situation
of (\ref{esplit3}), where $\fs_0'$ is the pairing in the lower
right corner.

\begin{equation}\label{esplit3}
\begin{split}
\epsfysize=3.8cm \epsffile{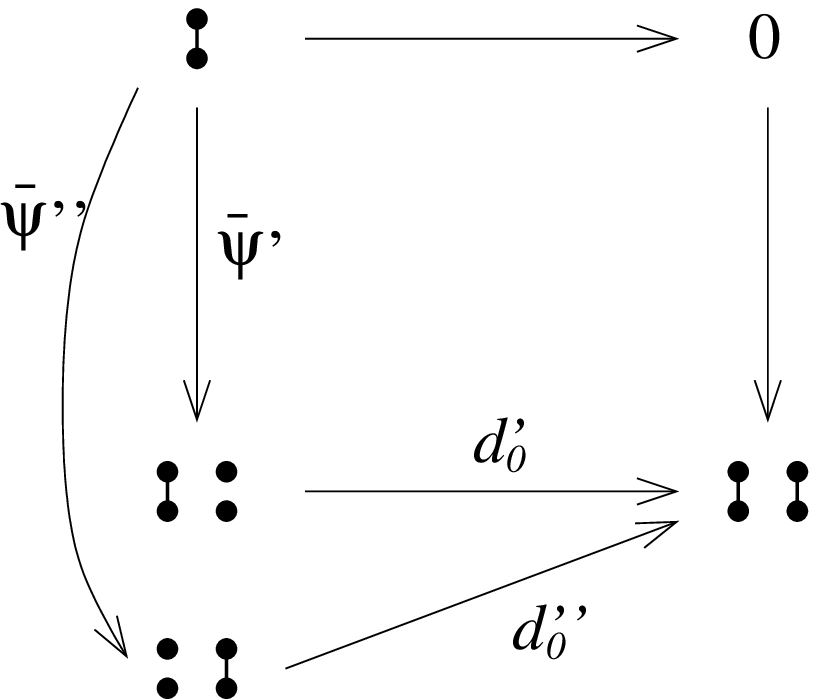}
\end{split}
\end{equation}
\v8

Now assume that $s_1'$ and $s_2'$ have no common
pair. Let $s_1$ and $s_2$ denote the intersections $s_1:=s\cap s_1'$
and $s_2:=s\cap s_2'$. Let $\fs_0$ denote the pairing of
the $\fn_0$--cable of $D_0$ which restricts to $s_1$ and $s_2$
on the components $K_1$ and $K_2$ and which agrees with $\fs$ on
all  other components of $D_0$. Then every pair of $\fs_0$ is also
be a pair of $\fs_0'$, and there is a unique pair $\pi$ of $\fs_0'$
which is not contained in $\fs_0$. We assume that $\pi$ belongs
$K_1$ or $K_2$ (for otherwise
$(d_0\bar{\psi})^{\fs_0',\fs} =\pm(\bar{\psi} d)^{\fs_0',\fs}$
is trivially satisfied).
If $\pi$ is disjoint from all pairs of $s_1\cup s_2$, we
are in the situation of (\ref{esplit1}), where $\pi$ is the pair
in the lower right corner.

\begin{equation}\label{esplit1}
\begin{split}
\epsfysize=3.4cm \epsffile{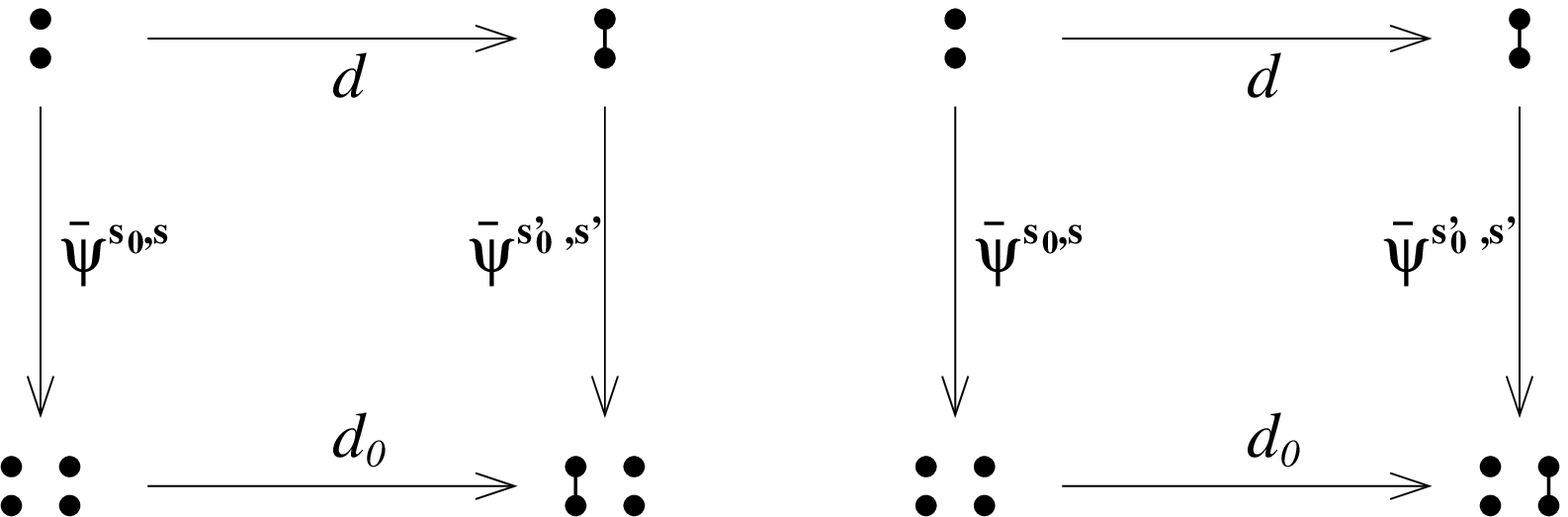}
\end{split}
\end{equation}
\v8

It is also possible that $\pi$ has a common dot with
a pair of $s_1\cup s_2$. Examples of this case are
shown in (\ref{esplit2}).

\begin{equation}\label{esplit2}
\begin{split}
\epsfysize=3.7cm \epsffile{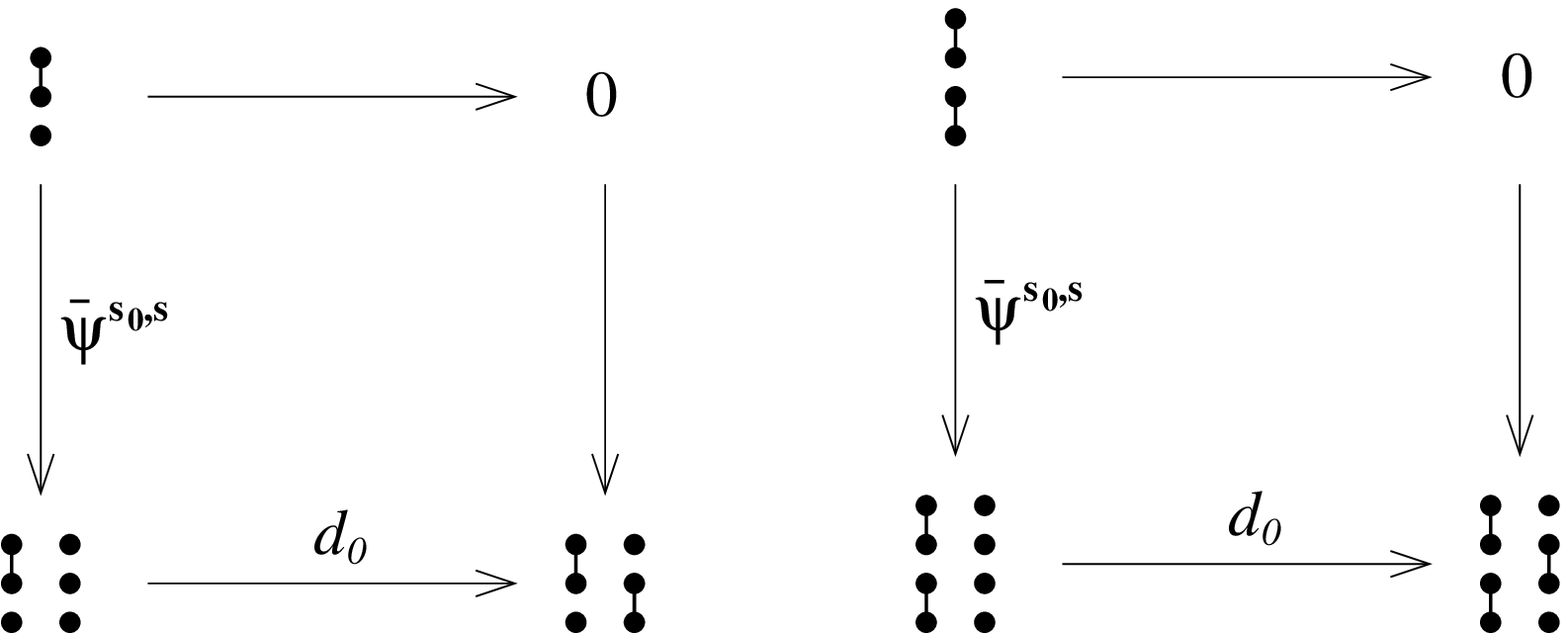}
\end{split}
\end{equation}
\v8

\begin{lem}\label{lssaddle}
Assume that the squares of (\ref{esplit1}) commute, up to sign,
and assume $d_0\bar{\psi}^{\fs_0,\fs}=0$ for all squares as in
(\ref{esplit2}). Then there is a
$0$--cochain $\gamma\in C^0(\G_{\fn_0}, \Z/2\Z)$ such that the morphisms
$(-1)^{\g(\fs_0)}\bar{\psi}^{\fs_0,\fs}$ determine a chain transformation
from $\ll D_\fn \rr_{\a,\b}$ to $\ll D_{0,\fn_0}\rr_{\a,\b}$.
\end{lem}
\begin{proof}
The proof is analogous to the proof of Lemma \ref{lmsaddle}
(although now we have to consider a map $f$ going from a subgraph $\G_{\fn_0}'$
of $\G_{\fn_0}$ to $\G_{\fn}$).
Note that the morphisms $d_0'\bar{\psi'}$ and $d_0''\bar{\psi'}$ of
(\ref{esplit3}) cancel automatically if the squares of (\ref{esplit1}) commute.
To see this, observe that the morphisms $\bar{\psi'}$ and $\bar{\psi''}$ of
(\ref{esplit3}) also appear in the squares of (\ref{esplit1}). Moreover, the
lower edges of the squares in (\ref{esplit3}) and (\ref{esplit1}) form a square
of $\G_{\fn_0}$. Now use that the squares of $\G_{\fn_0}$
anticommute.\end{proof}

\section{Towards functoriality}
Throughout this section let $\A$ be the category
of $\Z_{(2)}$--modules. 
%$\Kk(D_{\fn}):=\ll D_{\fn}\rr_{\a,\b}=\lr L_{\fn}\rl$.
Recall that $\F_{\rm Lee}, \F_{\Kh}:\Cobl^3\to \A$ extend to the
functors $\F_{\rm Lee}, \F_{\Kh}: \Kom(\Mat(\Kobh))\to \Kom(\Komh(\A))$.
The isomorphism classes of $\F_{\rm Lee}(\ll D_n\rr_{\a,\b})$ and
$\F_{\Kh}(\ll D_n\rr_{\a,\b})$
are invariants of a colored link.

\begin{thm}\label{tfunctor}
For $\a=\b=1$, the maps
$\F_{\rm Lee}(\psi^{\fs_0,\fs}_{1,1})$
and $\F_{\rm Lee}(\bar{\psi}^{\fs_0,\fs}_{1,1})$
induce chain transformations. 
\end{thm}

Before we prove the theorem, let us introduce a new
relation in $\Cob^3$, called the genus reduction
relation. Consider a cobordism $C'$
obtained from a cobordism $C$ by attaching two
small handles to a disk of $C$. The genus reduction
relation asserts that $C'=4C$. Now let us assume that
$2$ is invertible and that the relations (S), (T) and
(4Tu) hold. Then the genus reduction relation becomes
equivalent to the relation $\epsg{9mm}{g3}=8$, i.e. to the
defining relation for Lee's functor. As a consequence,
$\F_{\rm Lee}(\Kh(H(P)))^2/4$ is the identity map,
for $H(P)$ defined as in Subsection \ref{x}.

\begin{proof}[Proof of Theorem \ref{tfunctor}.]
Assume $\a=\b=1$ and assume that the genus reduction
relation holds. We have to show that under this assumption,
the morphisms $\psi^{\fs_0,\fs}_{1,1}$ and $\bar{\psi}^{\fs_0,\fs}_{1,1}$
satisfy the conditions of Lemmas \ref{lmsaddle} and \ref{lssaddle},
respectively. We start by proving 
$d_0\bar{\psi}_{1,1}^{\fs_0,\fs}=0$ for
the left square of (\ref{esplit2}) (the proof for the right square is
analogous). We assume that the three dots in the lower 
left corner of the square
are numbered from bottom to top from $i$ to $i+2$. Consider the diagram
$D^{\fn_0}_0$ of the $\fn_0$--cable of $D_0$. Let $l_i,l_{i+1}$,
and $l_{i+2}$ be three parallel edges of $D^{\fn_0}_0$, belonging
to the strands $C_i,C_{i+1}$ and $C_{i+2}$, respectively.
Let $G_i$ denote region of
$D^{\fn_0}_0$ which lies between $l_i$ and $l_{i+1}$, and let $G_{i+1}$ denote
the region which lies between $l_{i+1}$ and $l_{i+2}$. Choose a point $P_i$ on
$l_i$ and a point $P_{i+1}$ on $l_{i+1}$.
Observe that on $C_{i+2}$, the map $\bar{\psi}_{1,1}^{\fs_0,\fs}$
is induced by
a saddle cobordism. On $C_i$ and $C_{i+1}$, it is
induced by an annulus postcomposed with $(1+(\sigma(G_i)/2)\Kh(H(P_i)))$.
We can replace $\Kh(H(P_i))$ by $\Kh(H(P_{i+1}))$ because
we can move the point $P_i$ across the annulus.
For $\a=\b=1$, the map $d_0$ is given by
$(1+(\sigma(G_{i+1})/2)\Kh(H(P_{i+1})))$, postcomposed with an
annulus. Since $G_i$ and $G_{i+1}$ are neighbors,
$\sigma(G_{i+1})=-\sigma(G_i)$.
Summarizing, we see that $d_0\bar{\psi}_{1,1}^{\fs_0,\fs}$
factors through
$(1-(\sigma(G_i)/2)\Kh(H(P_{i+1})))
(1+(\sigma(G_i)/2)\Kh(H(P_{i+1})))
=1-\Kh(H(P_{i+1}))^2/4=0$.

To show that the squares of (\ref{emerge}) and (\ref{esplit1}) commute (up to
sign), apply isotopies, the neck cutting relation and the genus reduction
relation to the cobordisms corresponding to
$(d_0\psi_{1,1})^{\fs_0',\fs}$, $(\psi_{1,1} d)^{\fs_0',\fs}$,
$(d_0\bar{\psi}_{1,1})^{\fs_0',\fs}$ and $(\bar{\psi}_{1,1}d)^{\fs_0',\fs}$.
Use that $(\sigma(G)/2)\Kh(H(P))$ is independent of the choice of the
point $P$ on $C_i$, and that $\Kh(H(P))$ commutes with morphisms
induced by cobordisms which agree with the
identity cobordism in a neighborhood of $P$.
The details are left to the reader.
\end{proof}

\begin{thm}\label{tfu}
For $\a=0$, $\b=1$, the maps
$\F_{\Kh}(\psi^{\fs_0,\fs}_{0,1})$
and $\F_{\Kh}(\bar{\psi}^{\fs_0,\fs}_{0,1})$
induce chain transformations. 
\end{thm}

\begin{proof}
Let us first show that
 $\F_{\Kh}(d_0\bar{\psi}_{0,1}^{\fs_0,\fs})=0$ for
the left square of (\ref{esplit2}) (the proof for the right square is
analogous).
By the same arguments as in the previous proof,
 $\F_{\Kh}(d_0\bar{\psi}_{1,1}^{\fs_0,\fs})$
factors through
$\F_{\Kh}(\Kh(H(P_{i+1}))\Kh(H(P_{i+1})))
=0$, because the genus of the composition is bigger than one.

To show that the squares of (\ref{emerge}) and (\ref{esplit1}) commute (up to
sign), we have to proceed like in the previous proof, replacing
the genus reduction relation by the relation setting all cobordisms of genus
bigger than one to zero. An illustration in the case of the unknot
is given in Figure 5.
\end{proof}

\v8
\begin{center}
%\hspace*{-2.6cm}
\mbox{\epsfysize=5cm \epsffile{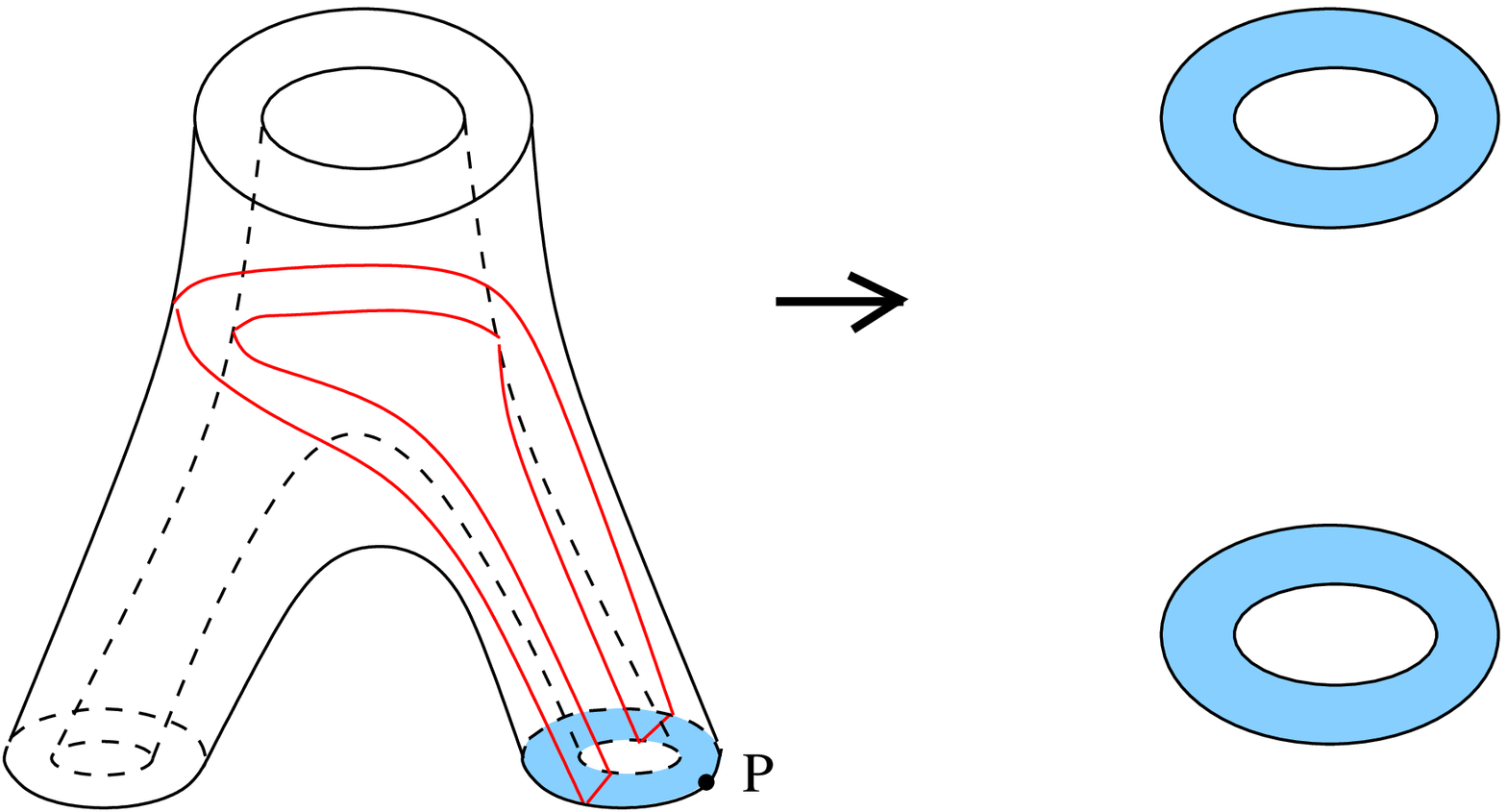}}
\v8

\parbox{12cm}{
Figure 5 {\it Commutativity of diagrams (\ref{esplit1}).} 
The filled regions are annuli of genus 1.
Applying the neck cutting relation to the red line,
and removing components of genus 2, we get
one half of the two genus 1 annuli, shown on the right. }
\end{center}

\begin{proof}[\bf Proof of Theorem 4]
>From Theorems \ref{tfunctor}, \ref{tfu} we know that saddles
induce chain transformations. It remains to show that the morphisms
associated to cups and caps do this also.
%induce chain transformations.
The case of cups is easy, so we only discuss the case of caps.
Let $D$ and $D_0$ be two link diagrams which are related by a cap
cobordism. Let $K$ be the component of $D_0$ which is not
contained in $D$. To show that
$\iota^{\fs_0,\fs}:\Kh(D^{\fs})\rightarrow\Kh(D_0^{\fs_0})$ induces a chain
transformation, we write the differential of
$\ll D_{0,\fn_0}\rr_{\a,\b}$ 
as a sum $d_0=d_0'+d_0''$, where $d_0'$ is the sum of
all morphisms 
 $\Kh(e)$ which increase the number of pairs on $K$, and $d_0''$
is the sum of all morphisms $\Kh(e)$ which increase the number of pairs on
one of the other components of $D_0$. It is easy to see that $d_0''$ commutes
with $\iota$. 
For $\a=0$, $\b=1$,
$\F_{\Kh}(d_0')\F_{\Kh}(\iota)=0$, because the genus of the
composition is bigger than one.
To complete the proof, we show that 
for $\a=1$ and $\b=1$, $\F_{\rm
Lee}(d_0')\F_{\rm Lee}(\iota)=0$. Note that $d_0'=\sum_{i=1}^n\Kh(e_i)$,
where $\Kh(e_i)$ is the morphism $(1 +X(K,i))$ composed with an
annulus glued to the strands $i$ and $i+1$ of the $n$--cable of
$K$. As in the proof of Theorem \ref{tfunctor}, 
we can replace $X(K,i)$ by $-X(K,i+1)$.
Recall that $\iota^{\fs_0,\fs}$ is given by
a union of cap cobordisms composed with the morphism
$\varphi=\sum_{j=0}^nA_j\circ B_j$. Using the genus reduction relation, we
obtain $\F_{\rm Lee}(( 1 +X(K,i))\circ A_j)=0$ for $i\leq j$ and $\F_{\rm Lee}((
1 -X(K,i+1)) \circ B_j)=0$ for $i\geq j$. Therefore $\F_{\rm Lee}(\Kh(e_i)\circ
A_j\circ B_j)=0$ for all $i,j$, and hence $\F_{\rm Lee}(d_0')\F_{\rm
Lee}(\iota)=0$.
\end{proof}

\begin{rem}
We do not know how to extend the original  colored
Khovanov bracket $\ll D_\fn\rr$
to a functor from the category $\CC^4$ to the category $\Kom(\Komh(\A))$.
%%%%%%%%%%%%%%%%%%%%%%%%%%%%%%%%%%%%%%%%%%%%%%%%%%%%%%%%%%%%%%%%%%%%%%%%%
For this bracket,
the morphisms $\psi^{\fs_0,\fs}_{0,2}$ induce chain transformations
(compare \cite{Kho}), but there is no choice of $\g,\delta$ for which the
morphisms $\bar{\psi}^{\fs_0,\fs}_{\g,\delta}$ induce chain transformations.
\end{rem}

\section{Rasmussen invariant for links}

\subsection{Definition}
By Lee's theorem \cite{Lee},
the homology of  $\F_{\rm Lee}(\ll L\rr)$ has rank $2^{|L|}$,
where $|L|$ is the number of components of $L$.
The generators of  Lee homology are in bijection
with the  orientations of $L$.
Hence, Lee homology of a knot has two generators.
In particular, the graded module associated to
the Lee homology of a knot
has two homogeneous generators, whose
topological degrees we denote
by $s_{\rm max}$ and $s_{\rm min}$.
The Rasmussen invariant $s(K)$ of a knot $K$ is
$$s(K):=\frac{s_{\rm max}+s_{\rm min}}{2}\, .$$

Let us extend the
 Rasmussen construction  to links.
Let   $L$ be an {\it oriented} link. Let $\s_o$ and $\s_{\bar o}$
 be
 the generators of the Lee homology corresponding
to the  orientation of $L$ and the opposite orientation,
respectively.  Then by Lemma 3.5 in \cite{Ra},
the filtered topological degrees of  $\s_o+\s_{\bar o}$
and $\s_o-\s_{\bar o}$ differ by two modulo 4.
Further, we  can show that they  differ by exactly  two.
Indeed,
a genus 1 cylinder  cobordism induces
an automorphism of $\F_{\rm Lee}(\ll L\rr)$ of topological degree $-2$,
which interchanges   $\s_o+\s_{\bar o}$ and $\s_o-\s_{\bar o}$.
The Rasmussen invariant $s(L)$ of a  link $L$
 is
$$s(L):=\frac{{\rm deg}(\s_o+\s_{\bar o})+{\rm deg}(\s_o-\s_{\bar o})}{2}.$$

%%%%%%%%%%%%%%%%%%%%%%%%%%%%%%%%%%%%%%%%%%%%%%%%%%%%%%%%%%%%%%%%%%%%%%%%%%%
 Note that $s(L)={\rm min}({\rm deg}(\s_o+\s_{\bar o}),
{\rm deg}(\s_o-\s_{\bar o}))+1$ and 
%%%%%%%%%%%%%%%%%%%%%%%%%%%%%%%%%%%%%%%%%%%%%%%%%%%%%%%%%%%%%%%%%%%%%%%%%%%
that the Rasmussen invariant of the $n$--component unlink
is $1-n$.

\subsection{Properties}
Let $S$ be a smooth oriented cobordism from $L_1$ to $L_2$ such that
every connected component of $S$ has boundary in $L_1$.
We will always  assume that the orientations of $L_1$ and $L_2$
coincide with ones induced by $S$, in the sense that $\partial S=-L_1\amalg
 L_2$.
Then the Rasmussen estimate generalizes to
\be\label{re}
|s(L_2)-s(L_1)|\leq -\chi(S)\ee
where $\chi(S)$ is the Euler characteristic of $S$. Indeed, arguing
as in \cite{Ra} we obtain the estimate $s(L_2)\geq s(L_1)+\chi(S)$.
By reflecting $S\subset \R^3\times [0,1]$ along $\R^3\times\{1/2\}$,
we obtain a cobordism from $L_2$ to $L_1$ with the same
Euler characteristic as $S$.
This gives us the estimate $s(L_1)\geq s(L_2)+\chi(S)$.

%%%%%%%%%%%%%%%%%%%%%%%%%%%%%%%%%%%%%%%%%%%%%%%%%
\vspace{0.5cm}

\begin{lem}\label{lproperties}
Let $\bar L$ be the mirror image of $L$ and $\#$, $\amalg$ denote the
connected sum and the disjoint union, respectively. Then
\begin{eqnarray}
\label{eunion} s(L_1\amalg L_2)&=&s(L_1)+s(L_2)-1 \\
\label{esum} s(L_1)+s(L_2)-2&\leq&s(L_1\#L_2)\;\:\leq\;\: s(L_1)+s(L_2)\\
 \label{emirror} -2|L|+2&\leq &s(L)+s(\bar L)\;\:\leq\;\: 2
%\label{esum} s(L_1\#L_2)&=&s(L_1)+s(L_2)         \\
% \label{emirror} -2|L|+2&\leq &s(L)+s(\bar L)\leq 0
\end{eqnarray}
\end{lem}

Note that the first inequality of (\ref{emirror}) becomes
an equality if $L$ is an unlink.
%The second inequality is an equality
%if $L$ is a Hopf link. Indeed, the Rasmussen
%invariant of the positive Hopf link is $+1$
%and the Rasmussen invariant of the negative Hopf link is
%$-1$.

\begin{proof}[Proof of Lemma \ref{lproperties}]
Let $o_1,o_2$ and $o$ denote the orientations of
$L_1,L_2$ and $L_1\amalg L_2$, respectively.
%%%%%%%%%%%%%%%%%%%%%%%%%%%%%%%%%%%%%%%%%%%%%%%%%%%%%%%%%%%%%%%%%%

The filtered modules $\F_{\rm Lee}(\ll L_1\amalg L_2\rr)$ and
$\F_{\rm Lee}(\ll L_1\rr)\otimes \F_{\rm Lee}(\ll L_2\rr)$
are isomorphic by an isomorphism
which sends $\s_o$ to $\s_{o_1}\otimes \s_{o_2}$.
Hence (\ref{eunion}) follows
from ${\rm deg}(\s_o)={\rm min}({\rm deg}(\s_o+\s_{\bar o}),
{\rm deg}(\s_o-\s_{\bar o}))=s(L_1\amalg L_2)-1$ and
${\rm deg}(\s_{o_i})={\rm min}({\rm deg}(\s_{o_i}+\s_{\bar {o_i}}),
{\rm deg}(\s_{o_i}-\s_{\bar {o_i}}))=s(L_i)-1$
(cf. \cite[Corollary 3.6]{Ra}).
(\ref{esum}) follows from (\ref{re}) and
(\ref{eunion}) because $L_1\amalg L_2$ and $L_1\# L_2$ are related by a
saddle cobordism. Similarly, (\ref{emirror}) can be deduced from (\ref{re})
and (\ref{eunion}) because there is a cobordism, consisting of $|L|$
saddle cobordisms, which connects $L\amalg\bar{L}$ to the
$|L|$--component unlink. 
%%%%%%%%%%%%%%%%%%%%%%%%%%%%%%%%%%%%%%%%%%%%%%%%%%%%%%%%%%%%%%%%%%
%Similarly, equation
%(\ref{esum}) follows
%because $\F_{\rm Lee}(\ll L_1\#L_2\rr)$ is isomorphic
%to $\F_{\rm Lee}(\ll L_1\rr)\otimes_V \F_{\rm Lee}(\ll L_2\rr)$.
%Here, $V$ denotes the ring $V:=\F_{\rm Lee}(\ll U\rr)\{-1\}$ where $U$ is the 
%unknot and the
%multiplication  is induced by the saddle
%cobordism. 
%We regard $\F_{\rm Lee}(\ll L_1\rr)$ and $\F_{\rm Lee}(\ll L_2\rr)$
%as $V$--modules by choosing distinguished
%points on $L_1$ and $L_2$ near the connected
%sum point (see \cite{Kh2} for details,
%where the notation $H^1$ is used
%instead of $V$, and where the topological grading is opposite to the
%the topological grading used in this paper).
%To show (\ref{emirror}), observe that there is a cobordism $S$
%from $L\#\bar{L}$ to  the $|L|$--component unlink, consisting
%of one cylinder and $|L|-1$ pair of pants surfaces.
%Combining (\ref{esum}) and (\ref{re}), we obtain
%$|s(L)+s(\bar{L})+|L|-1|=|s(L\#\bar{L})+|L|-1|\leq-\chi(S)=|L|-1$,
%which is equivalent to (\ref{emirror}).
\end{proof}

\subsection{Obstructions to 
 sliceness}
The notion of sliceness admits different generalizations to links.
We say that an oriented link $L$ is slice in 
{\it the weak sense} if there exists
an oriented smooth connected 
surface $P\subset B^4$ of genus zero, such that
$\partial P=L$.
$L$ is slice in {\it  the strong sense}
if every component bounds a disk in $B^4$ and all these
disks are disjoint.
Recently, Cimasoni and Florens \cite{CF} unified different
notions of sliceness by introducing colored  links.

The Rasmussen invariant of links is an obstruction 
to  sliceness.

\begin{lem}
Let $L$ be slice in the weak sense,
then
$$|s(L)|\leq  |L|-1 .$$
\end{lem}
\begin{proof}
If $L$ is slice in the weak sense,
then there exist an oriented genus 0 cobordism from $L$ to the unknot.
Applying (\ref{re}) to this cobordism we get the result.
\end{proof}

The multivariable Levine--Tristram signature defined in \cite{CF}
is also an obstruction to sliceness. However,
for knots with the trivial Alexander polynomial,
the Levine--Tristram signature is constant and equal to the
ordinary signature.
 Therefore, for a disjoint union of such knots 
the Rasmussen link invariant is often  a
 better obstruction 
 than the multivariable signature.
 Using the 
 Shumakovitch list of knots with the trivial Alexander polynomial,
but nontrivial Rasmussen invariant \cite{Sch}  and {\it  Knotscape},
one can easily construct
 examples. E.g. 
the multivariable signature of
$K_{15n_{28998}}\amalg K_{15n_{40132}}\amalg K_{13n_{1496}}$ vanishes
identically,
however
$s(K_{15n_{28998}}\amalg K_{15n_{40132}}\amalg K_{13n_{1496}})=4>3-1$, hence 
this split link
is not slice in the weak sense. 
Similarly, the Rasmussen invariant, but not the signature,
 is an obstruction to
sliceness for the following split links:
$K_{15n_{113775}}\amalg K_{14n_{7708}}$, $K_{15n_{58433}}\amalg
K_{15n_{58501}}$, etc.

\bibliographystyle{amsplain}

\begin{thebibliography}{10}
\bibitem{Ba} Baader, S.: {\it Braids and Rasmussen invariant},
  preprint (2006)
\bibitem{BN1} Bar--Natan, D.: {\it On Khovanov's categorification of
    the Jones polynomial}, Algebr. Geom. Topology {\bf 2} (2002) 337--370
\bibitem{BN} Bar--Natan, D.: {\it Khovanov's homology for tangles and
  cobordisms}, Geometry and Topology {\bf 9} (2005) 1443--1499
\bibitem{CS} Carter, S., Saito, M.: {\it Reidemeister moves for surface
isotopies and their interpretation as moves to movies}, J. Knot Theory
Ramific. {\bf 2} (1993) 251--284
\bibitem{CF} Cimasoni, D., Florens, V.: {\it Generalized Seifert surfaces
and signature of colored links}, to appear in Trans. AMS,
arXiv:math.GT/0505185
\bibitem{kh1} 
Khovanov, M.: {\it A categorification of the Jones polynomial},
Duke Math. J.
{\bf 101}  (2000)  359--426, arXiv:math.QA/9908171
%\bibitem{Kh2} Khovanov, M.:
% {\it A functor--valued invariant of tangles}, Algebr. Geom. Topology
%{\bf 2} (2002) 665--741, arXiv:math.QA/0103190
\bibitem{Kho} Khovanov, M.: {\it Categorifications of the colored
    Jones polynomial}, J. Knot Theory Ramific. {\bf 14} (2005) 111--130
\bibitem{Lee} Lee, E.: {\it On Khovanov invariant for alternating links},
arXiv:math.GT/0210213
%\bibitem{Li} Livingston, C.: {\it Computations of the Ozsv\'ath--Szab\'o
%knot concordance invariant}, Geometry $\&$ Topology {\bf 8} (2004) 735--742
\bibitem{MT} Mackaay, M., Turner, P.: {\it Bar--Natan's Khovanov
homology for coloured links}, arXiv:math.GT/0502445
\bibitem{Ra} Rasmussen, J.: {\it Khovanov homology and the slice
    genus}, arXiv:math.GT/0402131
\bibitem{Sch} Shumakovitch, A.: {\it Rasmussen invariant,
    slice--Bennequin inequality, and sliceness of knots},
arXiv:math.GT/0411643
\bibitem{W} Wehrli, S.: {\it Khovanov homology and Conway mutation}, 
arXiv:math.GT/0301312
\end{thebibliography}

\end{document}